\documentclass[letter,12pt]{article}
\usepackage{epsfig,psfrag,amsmath,amssymb,latexsym,eufrak,amscd,amsfonts,graphicx,fancyhdr}
\usepackage[margin=2cm,includeheadfoot]{geometry}
\pdfoutput=1

\newtheorem{theorem}{Theorem}[section]

\newtheorem{corollary}{Corollary}[section]

\newtheorem{example}{Example}[section]
\newtheorem{lemma}{Lemma}[section]
\newtheorem{proposition}{Proposition}[section]

\newcommand{\eop}{\mbox{ \vrule height7pt width7pt depth0pt}}

\numberwithin{equation}{section}

\title{Fluctuations of the longest common subsequence for sequences of independent blocks}
\author{Heinrich Matzinger\footnote{School of Mathematics, Georgia Institute of Technology, 686 Cherry Street, GA 30332-0160 Atlanta, USA}\quad Felipe Torres\footnote{Fakult\"at f\"ur Mathematik, Universit\"at Bielefeld, Postfach 100131 D-33501 Bielefeld, Germany}}
\date{\today}

\begin{document}
\maketitle
\thispagestyle{empty}
\abstract{\noindent 
The problem of the order of the fluctuation of the Longest Common Subsequence (LCS) of two independent sequences has been open for decades. There exist contradicting conjectures on the topic, \cite{CS} and \cite{W4}. Lember and Matzinger \cite{JJMM} showed that with i.i.d. binary strings, the standard deviation of the length of the LCS is asymptotically linear in the length of the strings, provided that $0$ and $1$ have very different probabilities. Nonetheless, with two i.i.d. sequences and a finite number of equiprobable symbols, the typical size of the fluctuation of the LCS remains unknown. In the present article, we determine the order of the fluctuation of the LCS for a special model of i.i.d. sequences made out of blocks. A {\it block} is a contiguous substring consisting only of one type of symbol. Our model allows only three possible block lengths, each been equiprobable picked up. For i.i.d. sequences with equiprobable symbols, the blocks are independent of each other. In order to study the fluctuation of the LCS in this model, we developed a method which reformulates the fluctuation problem as a (relatively) low dimensional optimization problem. We finally proved that for our model, the fluctuation of the ength of the LCS coincides with the Waterman's conjecture \cite{W4}. We belive that our method can be applied to any other case dealing with i.i.d. sequences, only that the optimization problem might be more complicated to formulate and to solve.
\section{Introduction}
\subsection{Motivation}\label{example}
In general trough this paper, $X$ and $Y$ will denoted two finite strings over a finite alphabet $\Sigma$. A common subsequence of $X$ and $Y$ is a subsequence which is a subsequence of $X$ as well as of $Y$. A Longest Common Subsequence of $X$ and $Y$ (denoted simply by {\rm LCS} of $X$ and $Y$, or only LCS when the context is clear enough) is a common subsequence of $X$ and $Y$ of maximal length. 
\\
\\
Let us motivate the study of the LCS of two string with an example: let $x=ACGTAGCA$ and $y=ACCGTATA$ two sequences over the finite alphabet $\Sigma=\{A,C,G,T\}$. A common subsequence of $x$ and $y$ could be $z=ATA$. For example, the string $z$ can be obtained from both $x$ and $y$ by just deleting some letters. We can represent the common subsequence $z$ as an {\it alignment with gaps} (a gap is denoted by '--'). The letters which are not in the subsequence get aligned with gaps, so that the subsequence has aligned the common letters of both sequences. The common subsequence $z=ATA$ can correspond to the following alignment:
\begin{equation}\label{alineamientohenry}
\begin{array}{c|c|c|c|c|c|c|c|c|c|c|c|c|c|c|c}
x&&A&C&-&-&G&-&T&A&G&-&T&-&A \\\hline
y&&A&-&C&C&-&G&T&A&-&C&-&A&-
\end{array}
\end{equation}
The representation of a subsequence as an alignment with gaps is not necessarily unique. However, each alignment with gaps defines exactly one common subsequence. We are interested only on alignments which aligns same-letter pairs or letters with gaps. In this paper, an alignment which aligns a maximum number of letter pairs of $x$ and $y$ is called {\it optimal alignment}. The subsequence defined by an optimal alignment is hence an LCS. The LCS of $x$ and $y$ is ${\rm LCS}(x,y)=ACGTAA$ and corresponds to the optimal alignment:
\begin{equation}\label{alineamiento3}
\begin{array}{c|c|c|c|c|c|c|c|c|c|c|c|c}
x&&A&C&-&G&T&A&G&-&T&A& \\\hline
y&&A&C&C&G&T&A&-&C&-&A&\\\hline
{\rm LCS}(x,y)&&A&C&&G&T&A&&&&A&
\end{array}
\end{equation}
In Bioinformatics (for instance \cite{W1,P1}), one of the main pro\-blems is to decide if two sequences  are related or not. If they are, it probably means that they evolved from a common ancestor. So, if they are related they should look somehow similar. Biologists try to determine which parts are related by finding an alignment which aligns the related parts. In our current example, the sequences $x=ACGTAGCA$ and $y=ACCGTATA$ are somehow similar, but if we compare them letter by letter the great similarity does not become obvious:
\begin{equation}\label{alineamiento1}
\begin{array}{c|c|c|c|c|c|c|c|c|c|c}
x&&A&C&G&T&A&G&T&A \\\hline
y&&A&C&C&G&T&A&C&A&
\end{array}
\end{equation}
In the alignment without gaps \ref{alineamiento1} we aligned mostly non-corresponding letter pairs, from where we obtained only $3$ aligned same-letter-pairs: (from left to right) the first, the second and the last pair. This is much less than what our optimal LCS alignment \ref{alineamiento3} achieved. A possible explanation why \ref{alineamiento1} is worse in looking for similarities than \ref{alineamiento3} is that some letters ``got lost'' in the evolution process, so that they are present only in one of the two sequences, so it is more useful to consider alignments with gaps instead to look for similarities. Longest Common Subsequences and Optimal Alignments \cite{W2,W3} are the main tools in computational biology to recognize when strings are similar. A relatively long LCS indicates that the strings are related, but how long does the LCS need to be to imply relatedness? Sequences which are not related are stochastically independent. Could it be that independent stochastic strings have a long LCS because of bad luck? To understand this questions, we need to figure out the size of the fluctuation of the LCS of independent strings. We are interested in the asymptotic of the fluctation since we mainly consider long sequences.
\subsection{Notation and history}
Let $X=X_1X_2\ldots X_n$ and $Y=Y_1Y_2\ldots Y_n$ be two stationary random sequences which are independent of each other, both drawn from the same finite alphabet $\Sigma$. Let $L_n:=|{\rm LCS}(X,Y)|$ denote the length of the LCS of $X$ and $Y$. A simple sub-additivity argument \cite{CS} shows that the expected length of the LCS divided by $n$ converges to a constant:
\[\lim_{n\rightarrow\infty}\frac{{\rm E}[L_n]}{n}=:\gamma>0.\]
The constant $\gamma$ depends on the distribution of $X$ and $Y$. But even for such simple cases as i.i.d. sequences with equiprobable symbols, the exact value of $\gamma$ is not known. Chv\`atal-Sankoff \cite{CS} derived upper and lower bounds for $\gamma$. These bounds were further refined by Baeza-Yates, Gavalda, Navarro and Scheihing \cite{BGNS}, Deken \cite{D}, Dancik-Paterson \cite{Paterson1,Paterson2} and finally Durringer, Hauser, Martinez, Matzinger \cite{martinezlcs,lcscurve}. The asymptotic value of the rescaling coefficient $\gamma$ as the number of symbols (the size of $\Sigma$) goes to infinity was determined by Kiwi, Loebl and Matousek \cite{KLM}. On the other hand, the speed of convergence was obtained by Alexander \cite{Alexander} by using techniques from percolation theory.
\\
The order of magnitude of the fluctuation of $L_n$ is unknown for situations as simple as i.i.d. sequences of equiprobable letters. In \cite{W4} Waterman conjectured that, in many situations, the fluctuation of the LCS is of order square root of the length times a constant:
\begin{equation}\label{ordern}
{\rm VAR}[L_n]=\Theta(n).
\end{equation}
Here the order $\Theta(n)$ means that there exist constants $0<a<b$ such that 
\[an\leq {\rm VAR}[L_n]\leq bn\]
for all $n\in\mathbb{N}$ (the constants $a$ and $b$ might depend on the distribution of $X$ and $Y$).
\\
\\
So far, Lember and Matzinger \cite{JJMM} proved the order given in \ref{ordern} for binary i.i.d. sequence, but when the probability of $1$ is much less than the probability of $0$. Durringer, Lember and Matzinger \cite{periodiclcs} obtained also the same order when one sequence is non-random, binary and periodic whilst the other binary sequence is i.i.d.  Bonetto and Matzinger \cite{BM} proved also the same order when the first sequence is drawn from a three letter alphabet $\{0,1,a\}$ whilst the second sequence is binary. Finally, Houdre and Matzinger \cite{HM} proved also the same orden when the two sequence are binary and i.i.d. but the {\it scoring function} which defines the alignment is such that one letter has a somewhat larger score than the other letter. Recall that in \cite{S}, Steele proved that there exists a constant $c>0$ not depending on $n$ such that ${\rm VAR}[L_n]\leq c\cdot n$, regardless of the alphabet $\Sigma$. This means that one only needs to find  good lower bounds for the variance of $L_n$ in order to find results on the fluctuation of $L_n$.
\\
\\
The LCS problem can be formulated as another popular open problem in probability theory, namely the Last Passage Percolation problem with correlated weights. The equivalence is as follows: let the set of vertices in our percolation setting be $V:=\{0,1,2,\ldots,n\}\times \{0,1,2,\ldots,n\}$. The set of oriented edges $E\subset V\times V$ contains horizontal, vertical and diagonal edges. The horizontal edges are  oriented to the right, whilst the vertical edges are oriented upwards. Both have unit length. The diagonal edges point up-right at a $45$-degree angle and have length $\sqrt{2}$. Hence $V:=\left\{\; (v,v+e_1),(v+e_2,v),(v,v+e_3) \;|\;v\in V\;\right\}$, where $e_1:=(1,0)$, $e_2:=(0,1)$ and $e_3:=(1,1)$. With the horizontal and vertical edges, we associate a weight of $0$. With the diagonal edge from $(i,j)$ to $(i+1,j+1)$ we associate the weight $1$ if $X_{i+1}=Y_{j+1}$ and $-\infty$ otherwise. In this manner, we obtain that the length of the LCS denoted by $L_n=|{\rm LCS}(X_1X_2\ldots X_n,Y_1Y_2\ldots Y_n)|$
is equal to the total weight of the heaviest path going from $(0,0)$ to $(n,n)$. Note that the weights on our 2-dimensional graph are not ``truly 2-dimensional'': they depend only on the one dimensional sequences $X=X_1\ldots X_n$ and $Y=Y_1\ldots Y_n$. 
\\
\\
The LCS problem is also related to the problem of the Longest Increasing Subsequence of a random permutation (for short only LIS), namely the LIS can be seen as the LCS of two sequences where one is a sequence of randomly permuted numbers and the other is the sequence of increasing integers. Take for example $5$ cards numerated from $1$ to $5$. Mix them thoroughly (until each permutation is equally likely). Then, lay them down face up in one line on a table. For example, you could obtain the permutation:
\[\begin{array}{ccccc}
2&3&1&5&4
\end{array}\]
A longest increasing subsequence here is $2,3,5$. We designate by $l_n$ the length of the longest increasing subsequence of such a random permutation, so in our case $l_5=3$. Note that the length of the LIS is equal to the length of the LCS of  the permutation and the sequence of increasing numbers. In our example $l_5=|{\rm LCS}(23154,12345)|.$ So, thanks to this relation and the recently tremendous breakthrough on the study of the LIS problem, many people was optimistic about finding a solution to the LCS problem by applying the new techniques from the LIS problem. Unfortunatelly, nobody has succeded so far in doing that. Moreover, we now belive that the LCS problem and the LIS problem are essentially from different classes though they have some features in common, for instance that both can be seen as passage percolation models, since the LIS problem is asymptotically equivalent to a special last passage percolation process on a Poisson graph. Let us recall some basic results about the LIS problem. In \cite{BaikDeiftJohansson99}, Baik, Deift and Johansson proved that
\[\frac{l_n-2\sqrt{n}}{n^{1/6}}\]
converges in distribution as $n\rightarrow\infty$ to a so called Tracy-Widom distribution (here $l_n$ denotes the length of the longest increasing sequence of a random permutation drawn from the symmetric group $\mathcal{S}_n$ with the uniform distribution). This limiting distribution can be obtained via the solution of  the Painleve II equation. It was first obtained by Tracy and Widom \cite{TW1,TW2} in the framework of Random Matrix Theory where it gives the limit distribution for the (centered and scaled) largest eigenvalues in the Gaussian Unitary Ensemble of Hermitian matrices. The problem of the asymptotic of $l_n$ was first raised by Ulam \cite{Ulam}. Substantial contributions to the solution of the problem have been made by Aldous and Diaconis \cite{Aldous99}, Hammersley \cite{Hammer}, Logan and Shepp \cite{Logan}, Vershik and Kerov (Vershik/Kerov 1977 Soviet math dokl). 
\\
\\
Coming back to the reason we belive make the LCS problem and the LIS problem essentially different, we can say the following: in the LIS case, the order of the fluctuation is power $1/3$ of the expectation and not square root. For the LCS case, the expectation is of order $n$. So, if the fluctuation was also a third power of the expectation, then we would have that ${\rm VAR}[L_n]$ should be of order linear in $n^{2/3}$. This is the order of magnitude conjectured by Chvatal-Sankoff \cite{CS} for which  several people  have some heuristic proofs. We believe that this order is wrong for the LCS-problem, based in all the above cites (and the present article) which confirmed Waterman's conjecture in many cases. However, for short sequences (small $n$) the order conjectured by Chvatal-Sankoff (which corresponds to the order of the fluctuation of the LIS) might be what one approximately observes in simulations. We believe that for short sequences, the underlying percolation structure shared between the LCS problem and the LIS problem make the two fluctuation look the same though the situation changes for large $n$: for short sequences, the correlation of the weights in the LCS problem has no strong effect and the system behaves as if the weights would be independent, as in the Poisson graph situation. So far, this arguments have not been rigoruosly proved, turning them in our opinion into attractive open questions in the area.

\section{Model and main ideas}\label{momi}


Let $l>0$ be an integer parameter. Let $B_{X1},B_{X2},\ldots$ and $B_{Y1},B_{Y2},\ldots$ be two i.i.d. sequences independent of each other such that:
\begin{eqnarray}
\mathrm{P}(B_{Xi}=l-1)=\mathrm{P}(B_{Xi}=l)=\mathrm{P}(B_{Xi}=l+1)&=&1/3 \nonumber\\
\mathrm{P}(B_{Yi}=l-1)=\mathrm{P}(B_{Yi}=l)=\mathrm{P}(B_{Yi}=l+1)&=&1/3. \nonumber
\end{eqnarray}
We call the runs of $0$'s and $1$'s blocks. Let $X^{\infty}=X_1X_2 X_3\ldots$ be the binary sequence so that the $i$-th block has length $B_{Xi}$ where $X_1$ is choosen $0$ with probability $1/2$ or $1$ with probability $1/2$. Similarly let $Y^{\infty}=Y_1Y_2 Y_3\ldots$ be the binary sequence so that the $i$-th block has length $B_{Yi}$ and $Y_1$ is choosen $0$ with probability $1/2$ or $1$ with probability $1/2$.

\begin{example}
Assume that $X_1=1$ and  $B_{X1}=2$, $B_{X2}=3$ and $B_{X3}=1$. Then we have  that the sequence  $X^\infty$ starts as follows $X^\infty=1100010\cdots$ meaning that in $X^\infty$ the first block consists of two 1's, the second block consists of three 0's, the third block consists of one 1's, etc.
\end{example}
Let $X$ denote the sequence obtained by only taking the first $n$ bits of $X^\infty$, namely $X=X_1X_2X_3\ldots X_n$ and similarly $Y=Y_1Y_2Y_3\ldots Y_n.$ Let $L_n$ denote the length of the {\rm LCS} of $X$ and $Y$, $L_n:=|{\rm LCS}(X,Y)|.$ 
\\
{\bf The main result of this paper states that for $l$ large enough, the order of the fluctuation  of $L_n$ is $n$}:
\begin{theorem} \label{maintheorem}
There exists $l_0$ so that for all $l\geq l_0$ we have that:
\[\mathrm{VAR}[L_n]=\Theta(n)\]
for $n$ large enough.
\end{theorem}
We show that the above theorem is equivalent to proving that ``a certain random modification has a biased effect on $L_n$''. This is a technique with similar approches in other papers (for instance see \cite{JJMM}, \cite{BM}). So the main difficulty is actually proving that the random modification has typically a biased effect on the {\rm LCS}. This random modification is performed as follows: we choose at random in $X$ a block of length $l-1$ and at random one block of length $l+1$, this means that all the blocks in $X$ of length $l-1$ have the same probability to be chosen and then we pick one of those blocks of length $l-1$ up and also that all the blocks in $X$ of length $l+1$ have the same probability to be chosen and we pick one of those blocks of length $l+1$ up. Then we change the length of both these blocks to $l$. The resulting new sequence is denoted by $\tilde{X}$. Let $\tilde{L}_n$ denote the length of the {\rm LCS} after our modification of $X$. Hence: 
\[\tilde{L}_n:=|{\rm LCS}(\tilde{X},Y)|.\]
If we can prove that our block length changing operation has typically a biased effect on the {\rm LCS} than the order of the fluctuation of $L_n$ is $\sqrt {n}$. This is the content of the next theorem:
\begin{theorem} \label{theorem2}
Assume that there exists $\epsilon>0$ and $\alpha>0$ not depending on $n$ such that for all $n$ large enough we have:
\begin{equation}\label{fundamental}
\mathrm{P}\left(\;\;  \mathrm{E}[\tilde{L}_n-L_n|X,Y]\geq \epsilon\;\;\right)\geq 1-\exp(-n^\alpha).
\end{equation}
Then,
$$\mathrm{VAR}[L_n]=\Theta(n)$$
for $n$ large enough.
\end{theorem}
The above theorem reduces the problem of the order of fluctuation to proving that our random modification has typically a higher probability to lead to an increase than to a decrease in score. The proof of this result is not included in the present article for shortness reasons, though all the details are in \cite{MT2}. In all what follows, we assume that theorem \ref{theorem2} is true. The next step is to ask: how can we prove, in our block model, that the condition \ref{fundamental} is satisfied? In theorem \ref{minimizingtheo}, we see that the condition \ref{fundamental} can be obtained from the positive solution of a minimizing problem. This minimizing problem has to do with the proportion of symbols which build up the LCS, been placed on a 9 dimensional space. By using Lagrange multiplyers techniques, we are able to further reduce it to a parametrized 3 dimensional optimization problem. Furthermore, we numerically and graphically verify that the positive minimum condition is already verified for $l>5$, which implies that ${\rm VAR}[L_n]=\Theta(n)$ holds already for $l=6$. Details on the solution of this minimization problem can be found in \cite{FelipePhD}.
\\
\\
The article is organized as follows: in what is left of section \ref{momi}, we explain how to relate the effect of the random modification with a constrained optimization problem on the prorportion of symbols used to build up the LCS and how this relation is used to prove theorem \ref{maintheorem}. In section \ref{combileftout}, we discuss some combinatorial aspects of aligned blocks in optimal alignments especific for this block model. Finally, in section \ref{details} we devote ourself to prove theorem \ref{minimizingtheo}

\subsection{Random modification and proportion of aligned blocks}
Let us next look, with the help of an example, when the random modification introduces an increase or a decrease in the score:
\begin{example}
Let us suppose $l=3$. Let us take two sequences $x=00110011110000111$ and $y=0011100001100001111$. An optimal alignment (in the sense of the Example \ref{example}) would be:
\begin{equation}
\label{alignment0}
\begin{array}{c|c|c|c|c| c|c|c|c|c| c|c|c|c|c| c|c|c|c|c| c|c|c}
 x&  &0&0&1&1&-&0&0&-&-&1&1&1&1&0&0&0&0&1&1&1&- \\\hline
 y&  &0&0&1&1&1&0&0&0&0&1&1&-&-&0&0&0&0&1&1&1&1 \\\hline
\mathrm{LCS}& &0&0&1&1& &0&0& & &1&1& & &0&0&0&0&1&1&1&
\end{array}
\end{equation}
In this example no block gets left out completely. By this we mean that no block is only aligned with gaps. The first block of $x$ is aligned with the first block of $y$. The second block of $x$ is aligned with the second block of $y$. By this we mean that all the bits from the second block of $x$  are either aligned with bits of the second block of $y$ or with gaps and vice versa. We have that  the second block of the {\rm LCS} is hence obtained from the second blocks of $x$ and $y$ by taking the minimum of their respective lengths. In our current special example, we have that for all $i=1,2,\ldots, 6$, the $i$-th block of $X$ gets aligned with the $i$-th block of $y$. We could represent this idea visually by viewing the alignment as an alignment of blocks in the following manner:
\begin{equation}
\label{alignment1}
\begin{array}
{c|c|l|l|l|l|l|l}
 x&  &00&11 &00 &1111&0000&111 \\\hline
 y&  &00&111&0000&11  &0000&1111 \\\hline
\mathrm{LCS}& &00&11 &00 &11  &0000&111
\end{array}
\end{equation}
Let us next analyze what is the expected change when we perform our random modification. In $x$ there are exactly $3$ blocks of length $l-1=2$. These are the first three blocks of $x$. The first block of $x$ of length $2$ is aligned with a block of $y$ of length $2$, the second one with a block of length $3$ and the fourth with a block of length $4$. Hence, when we increase the length of the first block of length $2$ of $x$ by one the score does not increase. When we increase the second or third, however, the score increases by one unit. Each of these blocks has the same probability $1/3$ to get drawn. Hence, the conditional expected  increase due to the enlargement of a randomly chosen block of length $2$ in this case, is equal to $2/3$. In our random modification we also choose a block of length $l+1$ and decrease it to length $l$. In our example, there are two blocks in $x$ of length $l+1=4$. These blocks are the fourth and fifth block of $x$. The fourth block is aligned with a block of length $2$ whilst the fifth is aligned with a block of length $4$. Hence, when we decrease the length of the fourth block we get no change in score whilst when we decrease the fifth we get a decrease by one unit. Each of the two blocks have same probability to get drawn. This implies that the expected change due to decreasing a randomly chosen block of length $4$ is equal to $-1/2$. Adding the two changes, we find that for $x$ and $y$ defined as in the current example, the conditional expected change is equal to:
\begin{equation}\label{simplified}
\mathrm{E}[\;\tilde{L}_n-L_n\;|X=x,Y=y]=\frac23-\frac12=\frac16
\end{equation}
In our example we have six aligned block pairs leading to the following set of pairs of lengths:
$$\{(2,2);(2,3);(2,4);(4,2);(4,4);(3,4)\}.$$
\end{example}
Let $p_{ij}$ designate the proportion of aligned block pairs which have the $x$-block having length $i$ and the $y$-block having length $j$. 
\begin{example}
For our example above we have:
\begin{equation}\left(
\begin{array}{ccc}
p_{22}&p_{23}&p_{24}\\
p_{32}&p_{33}&p_{34}\\
p_{42}&p_{43}&p_{44}
\end{array}
\right)=
\left(
\begin{array}{ccc}
\frac16&\frac16&\frac16\\
0      &0      &\frac16\\
\frac16&0      &\frac16
\end{array}
\right)
\end{equation}
\end{example}
\vspace{12pt}
With this notation, equality \ref{simplified}
can be written as:
\begin{equation}\label{simplified2}
\mathrm{E}[\;\tilde{L}_n-L_n\;|X=x,Y=y]\geq 
\frac{p_{l-1,l}+p_{l-1,l+1}}{p_{l-1,l-1}+p_{l-1,l}+p_{l-1,l+1}}-
\frac{p_{l+1,l+1}}{p_{l+1,l-1}+p_{l+1,l}+p_{l+1,l+1}}
\end{equation}
The inequality \ref{simplified2} holds if there exists an optimal alignment $a$ of $x$ and $y$ leaving out no blocks, and having a proportion $p_{ij}$ of aligned block pairs such that  the $x$-block has length $i$ and the $y$-block has length $j$ (for every $i,j\in\{l-1,l,l+1\}$). Typically, for large $n$, the optimal alignment will not be like in the example above, but there will be blocks which are left out, which implies also that some blocks are aligned with several blocks at the same time. Let us check an example:
\begin{example}
Let $x=00110011100011000$ and $y=00001111000011000$. In this situation the {\rm LCS} is equal to ${\rm LCS}(x,y)=000011100011000$ and corresponds to the following optimal alignment:
\begin{equation}
\label{alignment3}
\begin{array}
{c|c|c|c|c| c|c|c|c|c| c|c|c|c|c|c|c|c|c|c}
 x& 0&0&1&1&0&0&1&1&1&-&0&0&0&-&1&1&0&0&0  \\\hline
 y& 0&0&-&-&0&0&1&1&1&1&0&0&0&0&1&1&0&0&0  \\\hline
\mathrm{LCS}&0&0& & &0&0&1&1&1& &0&0&0& &1&1&0&0&0
\end{array}
\end{equation}
 which in block representation would be:
\begin{equation}
\begin{array}
{c|c|c|c|c|c}
 x& 001100&111 &000  &11&000 \\\hline
 y& 0000  &1111&0000 &11&000 \\\hline
\mathrm{LCS}&0000  &111 &000  &11&000
\end{array}
\end{equation}
In the last alignment above we see that the first block of $y$ is aligned with the first and third block of $x$. This implies that the second block of $x$ is ``completely left out'', which means all its bits are aligned with gaps. The other blocks are aligned one block with one block: the fourth block of $x$ is aligned with the second block of $y$, whilst the fifth block of $x$ is aligned with the third block of $y$. Finally the last blocks of $x$ and $y$ are aligned with each other. 
\end{example}
In everything that follows, the proportions $p_{ij}$ will only refer to the block pairs aligned one block with one block. Hence, in the alignment \ref{alignment3}, the first three blocks of $x$ and the first block of $y$ do not contribute to $\{p_{ij}\}_{i,j}$.

\begin{example}
In the last example above there are $4$ block-pairs aligned one block with one block. The corresponding pairs of block-lengths are:
$$(3,4);(3,4);(2,2);((3,3)$$
Hence for the alignment  \ref{alignment3}, we find $p_{3,4}=2/4$, $p_{2,2}=1/4$, $p_{3,3}=1/4$ and $p_{ij}=0$ for all $(i,j)\notin\{ (3,4),(2,2),(3,3)\}$. We will denote by $q_1$, resp. $q_2$, the proportion of left out blocks in $x$, resp. in $y$. In the alignment \ref{alignment3}, in the sequence $x$ there is one left out block from a total of $7$ blocks. This implies that $q_1=1/7$. There is no left out block in $y$ so that $q_2=0$. In section \ref{combileftout}, we will see that typically, for $n$ large enough, $q_1$ and $q_2$ can be taken as close to each other as we want to. When $q_1=q_2$ we denote the proportion of left out blocks by $q$. When we choose a block of length $l-1$ in $x$ to increase its length we will have to consider the probability that the block is not aligned one block with one block. In the alignment \ref{alignment3}, there are $4$ blocks in $x$ of length $l-1=2$. The first three are not aligned one block with one block: the second is left out, whilst the first and the third block are aligned with the same block of $y$. Hence in \ref{alignment4} the proportion of blocks not aligned one to one among the blocks of length $2$ is $3/4$. On the other hand, the blocks of length $l+1=4$ in $x$ are all aligned one to one. So, for the alignment \ref{alignment3}, we have that the proportion of blocks not aligned one to one among the blocks of length $4$ is $0$.
\end{example}
Using some combinatorial arguments, in section \ref{combileftout} we will see  that typically the proportion among the blocks of $x$ of length $l-1$ which are not aligned one block with one block is not more than $9q$. Similarly for the blocks of length $l+1$ in $x$ one gets a bound $3q$ for the proportion of  blocks aligned with several blocks of $y$ or left out. We can rewrite the lower bound on the right side of inequality \ref{simplified}, taking also into account the left out blocks. Assuming that there is an equal proportion of blocks $q$ which are not aligned one to one in $x$ and in $y$ we get the following lower bound for the conditional expected increase in the LCS:
\begin{equation}\label{lowerbound}
\frac{p_{l-1,l}+p_{l-1,l+1}}{p_{l-1,l-1}+p_{l-1,l}+p_{l-1,l+1}}(1-9q)-
\frac{p_{l+1,l+1}}{p_{l+1,l-1}+p_{l+1,l}+p_{l+1,l+1}}(1-3q)-3q
\end{equation}

\noindent The above lower bound for the conditional expected increase in LCS holds assuming that the following conditions holds:
\begin{itemize}
\item{}There exists an optimal alignment leaving out exactly the same proportion $q$ of blocks in $X$ and in $Y$. For that optimal alignment $a$, let $\{p_{ij}\}_{i,j}$ denote the empirical distribution of the aligned block pairs, so that $p_{ij}=P_{ij}(a)$.
\item{}There is exactly the same number of blocks in $X$ and in $Y$.
\item{}In $X$, each block lenght $l-1,l,l+1$ constitutes exactly $1/3$ of the blocks. Same thing in $Y$.
\end{itemize}
The above conditions do not typically hold exactly but only approximately. We first look at this somehow simplified case before looking at the general case (for the general case, see the proof of theorem 2.1.3). Let us next explain how we get the bound \ref{lowerbound} for this somehow simplified case (also, the reader should compare it to the version \ref{simplified2} with no gaps). Assume next that we have an optimal alignment $a$ with given empirical distribution $\{p_{ij}\}_{i,j\in\{l-1,l,l+1\}}$ of the aligned block pairs and leaving out in both sequences $x$ and $y$ a proportion $q$ of blocks.  What is now the effect of our random change on the score of the alignment $a$? First let us look at the randomly chosen block of length $l-1$ which gets its length changed to $l$. If that block is aligned with a block of length $l$ or $l+1$ the alignment gets increased by one unit. So, conditional that the randomly chosen block of length $l-1$ is a block aligned one to one, we get that the probability of an increase is equal at least to:
$$\frac{p_{l-1,l}+p_{l-1,l+1}}{p_{l-1,l-1}+p_{l-1,l}+p_{l-1,l+1}}.$$
Now, if the randomly chosen block of length $l-1$ is aligned with two or more blocks, then we also get an increase by one unit. If the chosen block however is aligned with a block of $Y$ which is aligned with several blocks of $Y$ (let us call it a {\it polygamist} block), then we have no increase. The same happens if the block is not aligned with a block of $Y$. There are at most a proportion of $3q$ blocks which are not aligned with any block or aligned together with polygamist block of $Y$.  There are about a proportion of $1/3$ blocks of length $l-1$. Hence among the blocks of length $l-1$, there is a proportion of at least $1-9q$ which are aligned one block with one block or aligned one with several. Hence we get that the conditional expected change due to changing the randomly chosen block of length $l-1$ to $l$ is equal at least to:
\begin{equation}\label{I}
\frac{p_{l-1,l}+p_{l-1,l+1}}{p_{l-1,l-1}+p_{l-1,l}+p_{l-1,l+1}}(1-9q).
\end{equation}
Similarly we can analyze the effect of the randomly chosen block of length $l+1$ which gets reduced to length $l$. If the block is aligned one block to one block and the length of the aligned block of $Y$ is $l+1$ then the score can get reduced by one. If the block is aligned with a block of $Y$ of length $l$ or $l-1$ the score does not get reduced. Hence, given that the block of length $l+1$ chosen is aligned one block to one block, the conditional expected change is not less than:
$$-\frac{p_{l+1,l+1}}{p_{l+1,l-1}+p_{l+1,l}+p_{l+1,l+1}}.$$
On the other hand, when the chosen block of length $l+1$ is aligned with several blocks of $Y$ then the score goes down by one unit. There are at most a proportion $q$ of blocks of $X$ aligned with several blocks of $Y$. So, among the blocks of length $l+1$ this represents a proportion of at most $3q$. Hence we get that at worst the expected change due to changing a random block from $l+1$ to $l$ is equal to:
\begin{equation} \label{II}
-\frac{p_{l+1,l+1}}{p_{l+1,l-1}+p_{l+1,l}+p_{l+1,l+1}}(1-3q)-3q
\end{equation}
Putting \ref{I} and \ref{II} together we get that the expected conditional change of the alignment score is bounded below as follows:

{\footnotesize \[\mathrm{E}[\Delta L_a|X,Y]\geq
\frac{p_{l-1,l}+p_{l-1,l+1}}{p_{l-1,l-1}+p_{l-1,l}+p_{l-1,l+1}}(1-9q)
-\frac{p_{l+1,l+1}}{p_{l+1,l-1}+p_{l+1,l}+p_{l+1,l+1}}(1-3q)-3q\]}
where $\Delta L_a$ denotes the change in score of the alignment $a$ due to the random modification of $X$.

\vspace{12pt}
\noindent Then, to prove inequality \ref{fundamental} in theorem \ref{theorem2}, it is thus sufficient to show that for all optimal alignments $a$ of $X$ and $Y$, expression \ref{lowerbound} is positive and bounded away from zero with high probability. Hence the next question is how can we prove that typically, for large $n$, expression \ref{lowerbound} is  larger than a positive constant not depending on $n$?
\begin{example}
Let us return back to the example of  alignment \ref{alignment3}. That alignment left out only one block, and that was the second block of $X$. We could now proceed in a different order. We could first decide which blocks get left out before generating the random sequences $X$ and $Y$. The resulting alignment is in general not optimal. On the other hand, such an alignment has the property that the block pairs aligned one to one are i.i.d. This is a very nice property for large deviation estimations, for instance. Let us give an example. Assume we request that the only left out block is the second block of $X$ (as in alignment \ref{alignment3}). Assume we redraw $X$ and $Y$ and obtain $X=001110011 000 11 000 $ and $Y=000 11110000111000$. Then we get as alignment and Common Subsequence (CS) the following:
\begin{equation}
\label{alignment4}
\begin{array}
{c|c|c|c|c| c|c|c|c|c| c|c|c|c|c| c|c|c|c|c| c|c}
 x& 0&0&1&1&1&0&0&1&1&-&-&0&0&0&-&1&1&-&0&0&0  \\\hline
 y& 0&0&-&-&-&0&-&1&1&1&1&0&0&0&0&1&1&1&0&0&0 \\\hline
\mathrm{CS}& 0&0& & & &0& &1&1& & &0&0&0& &1&1& &0&0&0
\end{array}
\end{equation}
which can be represented as an alignment of blocks by:
\[\begin{array}
{c|c|c|c|c|c}
 x& 0011100&11  &000 &11 &000  \\\hline
 y& 000    &1111&0000&111 &000 \\\hline
\mathrm{CS}& 000    &11  &000 &11  &000
\end{array}\]
\end{example}
In this case we use the term of common subsequence instead of the longest common subsequence because we are leaving a block out of the alignment, if we do not leave it out we might get a longer common subsequence (which does not happen in this case neither but might happens in the general case). So, in this last example, before drawing $X$ and $Y$, we know that the fourth block of $X$ gets aligned with the second block of $Y$ and this aligned pair builts the second block in the {\rm CS}. The length of the second block of the {\rm CS} has thus length equal to $\min\{B_{X4},B_{Y2}\}$. Similarly, before even drawing $X$ and $Y$, we know that the fifth block of $X$ gets aligned with the third block of $Y$. Hence, we have that the pair of lengths in the second  block pair is $(B_{X5},B_{Y3})$ whilst the third block of the CS has length $\min\{B_{X5},B_{Y3}\}$. Note that $(B_{X4},B_{Y2})$ is independent of $(B_{X5},B_{Y3})$ and $B_{X4}$ is independent of $B_{Y2}$ whilst $B_{X5}$ is independent of $B_{Y3}$. The distribution of each of the blocks $B_{X4}$, $B_{Y2}$, $B_{X5}$ and $B_{Y3}$ is unchanged, they take value $l-1$, $l$ or $l+1$ with equal probability $1/3$. Hence, $(B_{X4},B_{Y2})$ can take any of the nine values in the set $\{(i,j)|i,j=l-1,l,l+1\}$ with probability $1/9$.

\vspace{12pt}
\noindent When we specify an alignment by deciding which blocks we leave out before drawing $X$ and $Y$,  the aligned block pairs are ``almost'' i.i.d. Why do we say ``almost''? In the above example $(B_{X4},B_{Y2})$ and $(B_{X5},B_{Y3})$ are i.i.d. and not just close to be i.i.d. On the other hand,  block $B_{X7}$ in the case \ref{alignment3} is no longer in $X$ if the first, third and fourth blocks  get each increase by one unit. In this sense the blocks are not completely independent. But since we take $n$ large this is only a minor effect. We will take care of this detail in section \ref{details} and until then pretend that the aligned block pairs are i.i.d.

\vspace{12pt}
\noindent Note that for each alignment $a$ defined by specifying which blocks we left out before drawing $X$ and $Y$, the empirical distribution of the aligned blocks is random. We write $\{P_{ij}(a)\}_{i,j\in\{l-1,l,l+1\}}$ for this empirical distribution. Thus, $P_{ij}(a)$ denotes  the proportion of aligned block pairs where the block of $X$ has length $i$ and the block of $Y$ has length $j$. Given a non-random distribution $\{p_{ij}\}_{i,j\in\{l-1,l,l+1\}}$ we can ask what is the probability for the empirical distribution to be equal to the $\{p_{ij}\}_{i,j}$. The answer is, since the block pairs are close to i.i.d, the distribution is close to a multinomial distribution:
\begin{equation}\label{probpij}
\mathrm{P}\left( \;P_{ij}(a)=p_{ij},\forall i,j\in I_l  \; \right)
\approx \left( \begin{array}{c}
n^*\\
p_{l-1,l-1}n^*\;p_{l-1,l}n^*\;\ldots\;p_{l+1,l}n^*\;p_{l+1,l+1}n^*
\end{array}  
\right)  \left(\frac19\right)^{n^*}
\end{equation}
where $n^*$ designates the total number of aligned block pairs (here we act as if that number would be non-random). By using Stirling, the expression \ref{probpij} is approximately equal to:
\begin{equation}
\label{e^}
e^{(\ln(1/9)+H(p))n^*}
\end{equation}
where $H(p)$ designates the entropy of the empirical  distribution:
$$H(p)=\sum_{i,j\in\{l-1,l,l+1\}}p_{ij}\ln(1/p_{ij}).$$
A question arises: for a given aligned block pairs distribution $\{p_{ij}\}$, is it likely that there exist an alignment with that distribution and having a proportion $q$ of left out blocks? Let $\mathcal{A}(q)$ denote the set of alignments leaving out a proportion $q$ of blocks. Let $A$ denote the event that there exists an alignment in $\mathcal{A}(q)$ having its empirical distribution equal to $\{p_{ij}\}$. An upper bound for the probability $\mathrm{P}(A)$ is given  by the number of elements in $\mathcal{A}(q)$ times the probability \ref{probpij}. By using \ref{e^}, this product is close to:
\begin{equation}
\label{noname}
|\mathcal{A}(q)|\cdot e^{(\ln(1/9)+H(p))n^*}.
\end{equation}
But the size of the set $\mathcal{A}(q)$ is approximately equal to $e^{2H(q)n/l}$, since there are about $n/l$ blocks. Hence, expression \ref{noname} is approximately equal to:
\begin{equation}
\label{last}
 e^{(2H(q)n/l)+ (\ln(1/9)+H(p))n^*}. 
\end{equation} 
If we want the event $A$ to not have exponentially small probability in $n$, we need the logarithm of \ref{last} to be non-negative, which leads to the condition:
\begin{equation}\label{haupt}
2H(q)+(1-4q)(\ln(1/9)+H(p))\geq 0,
\end{equation}
where we used as lower bound on $n^*$ the number $(n/l)(1-4q)$.

\vspace{12pt}
\noindent We can now explain how we prove that typically, for all optimal alignment, expression \ref{lowerbound} is larger than a positive constant not depending on $n$. For this we simply need to find a $q_0$ so that we can prove that the optimal alignment leaves out at most a proportion of $q\leq q_0$ blocks and then show that expression \ref{lowerbound} is bounded away from zero under condition \ref{haupt} for  $q\in[0,q_0]$. 
\\
\\
\noindent Let $F^n(q)$ be the event that any optimal alignment of $X$ and $Y$ leaves out at most a proportion $q$ of bocks in $X$ and leaves out the same proportion $q$ of blocks in $Y$. In more details, given $q>0$ and an optimal alignment $a$ of $X$ and $Y$ in $F^n(q)$, we can count the number of blocks that are left out (not used in $a$) and divide this number by the total number of blocks in $X$ to obtain $q_1$, and also divide this number by the total number of blocks in $Y$ to obtain $q_2$, then we know that $q_1 \le q$ and $q_2 \le q$.
\begin{example}
Let us take again the case where $X=00111001100011000$ and \\ $Y=00011110000111000$, then we have as before the following common subsequence (CS) represented in an alignment:
\[\begin{array}
{c|c|c|c|c| c|c|c|c|c| c|c|c|c|c| c|c|c|c|c| c|c}
 X& 0&0&1&1&1&0&0&1&1&-&-&0&0&0&-&1&1&-&0&0&0  \\\hline
 Y& 0&0&-&-&-&0&-&1&1&1&1&0&0&0&0&1&1&1&0&0&0 \\\hline
\mathrm{CS}& 0&0& & & &0& &1&1& & &0&0&0& &1&1& &0&0&0
\end{array}\]
and represented as an alignment of blocks by:
\[\begin{array}
{c|c|c|c|c|c}
 X& 0011100&11  &000 &11 &000  \\\hline
 Y& 000    &1111&0000&111 &000 \\\hline
\mathrm{CS}& 000    &11  &000 &11  &000
\end{array}\]
Let us compute $q_1$ and $q_2$ in this case. For $X$ we have a total of 7 blocks and only 1 block is left out in the alignment, so $q_1 = 1/7$. For $Y$ we do not have left out blocks so $q_2=0$. Then given $q>0$, this alignment belongs to $F^n(q)$ if and only if $q_1=1/7 \le q$ and $q_2=0 \le q$.
\end{example}
The next theorem says that if we can bound expression \ref{lowerbound} away from zero under condition \ref{haupt}, then we have typically the desired bias for $\mathrm{E}[\tilde{L}_n-L_n|X,Y]$ the conditional expected increase in score: 

\begin{theorem} \label{minimizingtheo}
Assume that there exists $q_0\in[0,(1/3)[$ such that
the following minimizing problem: 
\begin{equation} \label{expressionI}
\min \left(\frac{p_{l-1,l}+p_{l-1,l+1}}{p_{l-1,l-1}+p_{l-1,l}+p_{l-1,l+1}}(1-9q)- \frac{p_{l+1,l+1}}{p_{l+1,l-1}+p_{l+1,l}+p_{l+1,l+1}}(1-3q)-3q\right)
\end{equation}
under the conditions: 
\begin{equation} \label{conditionq}
q\in[0,q_0],\sum_{j}p_{l-1,j}\geq((1/3)-q_0)/2\;,\; \sum_{j}p_{l+1,j}\geq ((1/3)-q_0)/2
\end{equation}
\begin{equation} \label{conditionpij1}
\sum_{i,j\in I}p_{ij}=1,p_{ij}\geq 0,\forall i,j\in I
\end{equation}
\begin{equation} \label{entropycondition}
2H(q)+(1-4q)\left( \ln(1/9)+H(p)\right) \geq 0
\end{equation}
has a strictly positive solution. Let this minimum be equal
to $2\epsilon>0$.
Then we have that:
\begin{equation} \label{enbeta}
\mathrm{P}\left(\;\;  \mathrm{E}[\tilde{L}_n-L_n|X,Y]\geq \epsilon\;\;\right)\geq
1-e^{-n^\beta}-\mathrm{P}(F^{nc}(q_0)) 
\end{equation}
where $\beta>0$ is a constant not depending on $n$.
\end{theorem}
%
Note that the high probability of the biased effect is only given when ${\rm P}(F^{nc}(q_0))$ is small (recall that $F^{n}(q_0)$ is the event that in any optimal alignment the proportion of left out blocks is less/equal to $q_0$). This means that, in order to apply the above theorem, we first need to come up with a way to bound the proportion $q$ of left out blocks in any optimal alignment. If the proportion of left-out blocks $q$ is too high, the joint distribution $p_{ij}(a)$ of the aligned block lengths could just be anything. In other words, the entropy condition \ref{entropycondition} becomes useless when $q_0$ is not small enough. We can now summarize how to apply the last theorem above: we first need to establish that the proportion of left out blocks is small enough. This means that we need to find a $q_0$ which satisfies that ${\rm P}(F^n(q_0))$ is close to $1$ and small enough, so that the objective function \ref{expressionI} is bounded away from $0$ under the constrains \ref{conditionq}, \ref{conditionpij1} and \ref{entropycondition}. In section \ref{leftoutblocks}, we show that the proportion of left out blocks does typically not exceed any $q_0$ for which $q_0>(4/9)/(l-1)$, where $l$ is the average block length. With this bound on $q$, (that is taking $q_0=(4/9)/(l-1)$) we are then able to verify numerically \ref{FelipePhD} that the objective function \ref{expressionI} is bounded away from $0$ under our constrains, already for $l=6$. This then implies that for any $l\geq 6$, the order of the fluctuation is $\textrm{VAR}[L_n]=\Theta(n)$.
In the next section, we prove the last theorem precisely, taking care of other details, for example that the proportion of left out blocks in $X$ and in $Y$ does not coincide in every alignment, only in the optimal alignment. Other important point is that the probability $\mathrm{P}(F^{nc}(q))$ depends on the 
parameter $l$. In chapter \ref{leftoutblocks}, we show how to find upper bounds on the proportion of left out blocks. In general, for $l$ larger, the bounds gets better. Actually the bounds even converge to zero as $l$ goes to infinity. As $q$ goes to zero, expression \ref{expressionI} gets close to $1/3$ on the domain. That is why the minimizing problem in theorem \ref{minimizingtheo} has a strictly positive solution when $l$ is large enough.
\\
\\
Let us no prove that theorem \ref{minimizingtheo} and theorem \ref{theorem2} together imply theorem \ref{maintheorem}:

\vspace{12pt}
\noindent {\bf Proof.} We suppose that  $F^{nc}(q_0)$ has exponentially small probability for any fixed $q_0>0$ provided $l$ is large enough (see section \ref{leftoutblocks} and \ref{nLOBOA}). In section \ref{details} we will show how large $l$ should be depending on $q_0$ but not on $n$. The conditions in theorem \ref{minimizingtheo} are satisfied when $q_0>0$ (hence $q\leq q_0$ small enough) is taken small enough. Let us explain why. First note that inequality \ref{entropycondition} can be written:
\begin{equation}\label{Hp}
H(p)\geq \frac{-2H(q)}{1-4q}+\ln(1/9)
\end{equation}
When $q$ goes to zero, then $H(q)$ also goes to zero and so does $2H(q)/1-4q$. But we have that $H(p)$ is always less or equal to  $\ln(1/9)$, with equality iff all the $p_{ij}$'s are equal to $1/9$.
\\
It follows that by taking $q>0$ small enough, we get condition \ref{Hp} to imply that the distribution $p_{ij}$ gets as close as we want to the equiprobable distribution. On the other hand, when $q$ goes to zero and all the $p_{ij}$'s converge to $1/3$, then the quantity
$$\frac{p_{l-1,l}+p_{l-1,l+1}}{p_{l-1,l-1}+p_{l-1,l}+p_{l-1,l+1}}(1-9q)-
\frac{p_{l+1,l+1}}{p_{l+1,l-1}+p_{l+1,l}+p_{l+1,l+1}}(1-3q)-3q,$$ 
converges to $2/3-1/3=1/3>0$. This shows that by taking $q_0>0$ small enough we get that the minimizing problem in theorem \ref{minimizingtheo}, has a strictly positive solution. So, assume that $q_0>0$ is such that the following two things hold:
\begin{itemize}
\item $F^{nc}(q_0)$ has exponentially small probability in $n$.
\item The minimizing problem in theorem \ref{minimizingtheo}, has a strictly positive solution. Call this solution $2\epsilon$, where $\epsilon>0$.
\end{itemize}
By theorem \ref{minimizingtheo}, we then have that inequality \ref{enbeta} holds. But since $F^{nc}(q_0)$ is exponentially small in $n$, we get that the expression on the right hand side of inequality \ref{enbeta} is smaller than $\exp(-n^{\alpha})$ for all $n$ and $\alpha>0$ not depending on $n$. This implies that condition \ref{fundamental} in theorem \ref{theorem2} is satisfied. Then theorem \ref{theorem2} implies that:
$$\mathrm{VAR}[L_n]=\Theta(n).\qquad\eop$$
\section{Combinatorics of the left out blocks} \label{combileftout}
In this section we describe some of the combinatorial properties that our block model exhibits when looking for optimal alignments. Let us begin with an example:
\begin{example}
Let $X=0011100$ and $Y=0001100$. The {\rm LCS} is $001100$.  This corresponds to the following alignment:
\begin{equation}
\label{alignmentI}\begin{array}{c|c|c|c|c|c|c|c|c|c}
X&  &0&0&-&1&1&1&0&0\\\hline
Y&  &0&0&0&1&1&-&0&0\\\hline
\mathrm{LCS}&&0&0& &1&1& &0&0
\end{array}
\end{equation}

\noindent In this example, the first block of the {\rm LCS} has length $2$. It is obtained from the first block of $X$  and the first block of $Y$. The first block of $X$ has length $2$ whilst the first block of $Y$ has length $3$. The length of the first block of the {\rm LCS} is equal to the minimum of these two numbers. In this kind of situation we say that the first block of $X$ is aligned to the first block of $Y$. Similarly the length of the second block of the {\rm LCS} is the minimum of the lengths of the second block of $X$ and of $Y$.  We say that in this alignment the second block of $X$ gets \textit{aligned} with the second block of $Y$. Finally the third block of $X$ gets aligned with the third block of $Y$ to yield the third block of the {\rm LCS}. In this present example no block of $X$ or $Y$ got left out completely: every block ``contributed'' some bits to the {\rm LCS}. All the blocks are aligned one block of $X$ with one block of $Y$. Each such pair of aligned blocks is responsible for one block in the {\rm LCS}.
\end{example}
In some other cases, some blocks of $X$ and $Y$ are completely left out. Let us look at such a situation. 
\begin{example}
Consider $X=00100000111$ and $Y=00000100011$. The {\rm LCS} would be $000000011$. The {\rm LCS} corresponds to the alignment:
\begin{equation}
\label{alignmentII}
\begin{array}{c|c|c|c|c|c|c|c|c|c|c|c|c|c|c}
X&  &0&0&1&0&0&0&-&0&0&-&1&1&1 \\\hline
Y&  &0&0&-&0&0&0&1&0&0&0&1&1&- \\\hline
{\rm LCS}&&0&0& &0&0&0& &0&0& &1&1&
\end{array}
\end{equation}
\end{example}
In the last example above we have that the second block of $X$ and of $Y$ are totally left out and do not contribute to the {\rm LCS}.  We say that these blocks are {\it left out} blocks. The last block of $X$ and the last block of $Y$ ``get aligned'' together to yield the last block of the {\rm LCS}. We say that this is an \textit{aligned block pair} or also that these two blocks are \textit{aligned one block to one block}. One way of thinking about the {\rm LCS} defined by the alignment \ref{alignmentII} above is as follows: we first decide which blocks we leave out in $X$ and $Y$. Then from the two obtained sequences, we align block by block without leaving out any blocks. So the alignment \ref{alignmentII} can be seen as the alignment in which we leave out the second block of $X$ and the second block of $Y$. This gives then the modified sequences $X^*=0000000111$ and $Y^*=0000000000011$. Then we align $X^*$ and $Y^*$ block by block. The common subsequence we obtain has its $i$-th block having length equal to the minimum of the length of block $i$ of $X^*$ and of $Y^*$. In this example we have that the first and the third block of $X$ get aligned with the first and third block of $Y$. By this we mean that in both sequences the first and third block are made into one block and these blocks are then matched. We will be able to prove that in the case we study here this is untypical: for optimal alignment we will only have one block aligned with several at the same time, but not several with several. Let us look at one more example:
\begin{example}
Let $X=001001111$ and $Y=000011011$. The {\rm LCS} is $00001111$. This corresponds
to the alignment:
\begin{equation}
\label{alignmentIII}
\begin{array}{c|c|c|c|c|c|c|c|c|c|c|c}
X&  &0&0&1&0&0&1&1&-&1&1\\\hline
Y&  &0&0&-&0&0&1&1&0&1&1\\\hline
{\rm LCS}&&0&0& &0&0&1&1& &1&1
\end{array}
\end{equation}
Here the second block of $X$ is left out. Hence the first and the third
block of $X$ get aligned with the first block of $Y$. Similarly
the fourth block of $X$ gets aligned with the second and fourth
block of $Y$. The third block of $Y$ is left out. 
\end{example}
This situation will
happen in optimal alignment: one block aligned with several blocks of the other sequence.

\vspace{12pt}
\noindent Assume that we know for an alignment $a$ which blocks are left out. Assume that $X^*$, resp. $Y^*$ denotes the modified sequence $X$, resp. $Y$ where we left out the specified blocks. Let $Z$ denote a common subsequence defined by the alignment $a$. The alignment must then align all the blocks of $X^*$ with the blocks of $Y^*$ one to one, otherwise there would be more left-out blocks. Hence, the first block of $X^*$ gets aligned, then the second block of $X^*$ and so on. If the alignment wants to stand a chance to be an optimal one (and hence $Z$  to be a {\rm LCS}) for each pair of aligned blocks from $X^*$ and $Y^*$ aligned to one another, it needs to extract a maximum of bits of each such pair. Hence, for every $i=1,2,\ldots,j$ we have that the length of the block number $i$ of $Z$ must be equal to the minimum between the length of the $i$-th block of $X^*$ and the length of the $i$-th block of $Y^*$ (here $j$ denotes the number of blocks in $Z$.) Hence, since we are interested in {\rm LCS}'s (and hence in optimal alignments) we will only consider alignments defined in the following manner: first we define exactly which blocks get left out. Second we align the resulting sequences $X^*$ and $Y^*$ one block with one block. The next lemma says that in our setting  an optimal alignment cannot align several blocks with several blocks.

\vspace{12pt}
\noindent Another useful fact is that for optimal alignments we do not need to consider adjacent left-out blocks  except maybe at the end of the sequences. But in section \ref{details} we prove that only a small percentage of bits could be left out at the end of $X$ and $Y$ in an optimal alignment. Hence, the practical implication is that we only need to consider left out blocks at least separated by one non-left out block. Let us first explain what we mean by adjacent left out blocks between aligned blocks:
\begin{example}
Take $x=11001100$ and $y=00001100$. Let us align as follows:
\begin{equation}
\label{arrayI}
\begin{array}{c|c|c|c|c| c|c|c|c|c| c|c|c|c}
x& &-&-&-&-&1&1&0&0&1&1&0&0\\\hline
y& &0&0&0&0&1&1&-&-&-&-&0&0 
\end{array}
\end{equation}
\end{example}
We see a typical situation where the second and third block of  $x$ in the alignment above get left out (i.e. entirely aligned with gaps). These two blocks are adjacent and they are comprised between aligned blocks (i.e. in our example they are comprised between the first and fourth block of $x$ which are ``aligned'', by aligned we mean aligned with another block hence not entirely aligned with gaps). The next lemma below states that for our {\rm LCS} problem (i.e. optimal alignments) the kind of situation we face in the current numerical example \ref{arrayI} can be discarded. The reason is as follows. In the current example $y_7$ gets aligned with $x_7$. Now instead align $y_7$ with $x_3$ and keep all the rest of alignment \ref{arrayI} identical otherwise. Then by doing this you have not decreased the score but have destroyed the situation of two adjacent completely left out blocks. The next lemma shows what we explained in our example in a rigorous way:

\begin{lemma}\label{AadjacentBlocks} 
There exists an optimal alignment of $X$ and $Y$ having no adjacent left-out blocks between aligned blocks.
\end{lemma}

\vspace{12pt}
\noindent {\bf Proof.} View an alignment as a finite sequence of points
in $\mathbb{N}\times\mathbb{N}$, so that if $x_i$ gets aligned with $y_j$,
then $(i,j)$ is a  point in the set representing the alignment.
Introduce for two alignments 
$a,b\in\mathbb{N}\times\mathbb{N}$ the order relation
$a\leq b$ iff all $a$ contains the same number of points
as $b$ and if we numerate in both sets the points from 
down left to up right then the $i$-th point $a_i=(a_{ix},a_{iy})$ of $a$ and the $i$-th point $b_i=(b_{ix},b_{iy})$ of $b$ satisfy
$a_{ix}\leq b_{ix}$ and $a_{iy}\leq b_{iy}$ for all 
$i\leq |a|$. Here $|a|$ designates the number
of points in $a$.
Take now an optimal alignment which is minimal according to
the relation $\leq$. That optimal alignment satisfies
the property of not having several adjacent left out blocks
between aligned blocks. $\quad\eop$

\vspace{12pt}
\noindent Next we show the relation between left out blocks at the end of each sequence and the total left out blocks in each sequence:
\begin{lemma}\label{delta1delta2}
Let $x,y\in\{0,1\}^n$ be two sequences of length 
$n$. Let the number of blocks of $x$, resp. $y$ be denoted
by $n^*_1=(n/l)+\Delta_1$, resp. $n^*_2=(n/l)+\Delta_2$. Assume that
$|\Delta_1|,|\Delta_2|\leq \Delta$. Assume also that
$a$ is an alignment of $x$ and $y$ which does never leave out
adjacent blocks except maybe a contiguous group at the very end of $x$ and of $y$. Let
$\delta_1\geq0$ denote the proportion of blocks
which are entirely left out at the end of $x$, resp. $y$, among all the blocks
of $x$, resp. $y$.
 Let $q_1$, resp. $q_2$ denote
the proportion of blocks left out in $x$, resp. in $y$.
Then we find that:
\begin{equation}
\label{blockinequ}
|q_1-q_2|\leq 1.5|\delta_1-\delta_2|+\frac{4l\Delta}{n}
\end{equation}
\end{lemma}

\vspace{12pt}
\noindent {\bf Proof.} Let $x^*$, resp. $y^*$ denote the sequence we obtain
after we removed the blocks which are completely left
out by $a$. Since there are no other completely left
out blocks, we have that the number of blocks in $x^*$
must be equal to the number of blocks in $y^*$.
Note that for every left out blocks which has no
adjacent left out block the number of blocks is reduced
by $2$. for the adjacent left out blocks at the end,
for each left out block there is one block less.
Since there are no adjacent left out blocks except
the adjacent blocks at the end, we get that the number of 
blocks of $x^*$, resp. of $y^*$ is equal to
\begin{equation}
\label{n1}n_1^*(1-2(q_1-\delta_1)+\delta_1)),
\end{equation}
resp.
\begin{equation}
\label{n2}n_2^*(1-2(q_2-\delta_2)+\delta_2)).
\end{equation}
Taking the difference of \ref{n1} and \ref{n2} and
dividing by $(l/2n)$, we find
\begin{equation}
\label{q1q2}
q_1-q_2=1.5(\delta2-\delta_1)+\frac{bl\Delta}{n}
\end{equation}
where 
$$2b=1-2(q_1-\delta_1)+\delta_1)-
\left(1-2(q_2-\delta_2)+\delta_2)\right)$$
we see that $b$ is always smaller than $4$
which ends the proof. $\quad\eop$
\begin{lemma}
For $l>4$ any optimal alignment of $X$ and $Y$ does not align several blocks in $X$ with several blocks in $Y$.
\end{lemma}

\vspace{12pt}
\noindent {\bf Proof.}
Let us explain the idea behind through an example. Let us take  $x=0001111000111100000$ and $y=0001111000001110000$ two realizations of $X$ and $Y$, respectively, with $l=4$. An alignment using all blocks of $x$ and $y$ in block representation becomes:
\begin{equation}\label{AN1}
\begin{array}{c|c|c|c|c|c|c}
x&&000&1111&000&1111&00000\\\hline
y&&000&1111&00000&111&0000\\\hline
{\rm LCS}&&000&1111&000&111&0000
\end{array}
\end{equation}
Let us now suppose that we leave out the second block of $x$ and the second block of $y$, then the alignment in block representation looks like:
\begin{equation}\label{AN2}
\begin{array}{c|c|c|c|c}
x&&0001111000&1111&00000\\\hline
y&&000111100000&111&0000\\\hline
{\rm LCS}&&000000&111&0000
\end{array}
\end{equation}
One clearly sees that in alignment \ref{AN2} we lost the entire block of 1's of length 4 and we did not gain any new aligned symbol, so the {\rm LCS} decreased on $4$ units compared to alignment \ref{AN1}. In this particular example, the neighbour blocks of the left out block in $y$ had all together at least as many symbols (8 zeros all together) as the neighbour blocks of the left out block in $x$ had all together (6 zeros). In general we could gain at most 2 new symbols from the neighbour blocks of the left out block but we always loose at least $l-1$ symbols leaving a block out and aligning its neighbours together instead. The other blocks do not get involved in the change on the score. Hence, when one leaves out a block and tries to align the neighbour blocks together the {\rm LCS} changes in $2-(l-1)=3-l$. Then for blocks of length $l>4$, to align several blocks with several blocks decreases the {\rm LCS} rather than to increase it.
$\quad\eop$
\subsection{Maximum number of left out blocks}\label{leftoutblocks}
The first key question is the percentage of blocks which are at
most left out in an optimal alignment.  Since the blocks have length $l-1$, $l$ or $l+1$ with
equal probability $1/3$ the expected block length is $l$. Hence,
the expected number of blocks in a sequence of length $n$ is
about $n/l$. Now let us define the limit:
\begin{equation}\label{gammal}
\gamma_l=\lim_{n\rightarrow\infty}\frac{\textrm{E}[L_n]}{n}.
\end{equation}
Hence, the number of bits in the sequence $X$ (and also in the sequence
$Y$) which are not used for the {\rm LCS} is about $(1-\gamma_l)n$. Every block
we leave out means at least $l-1$ non-used bits.  Hence,
the number of left out blocks for long sequences
can typically not be much above:
$$\frac{(1-\gamma_l)n}{l-1}.$$
This represents typically a proportion of:
$$\frac{(1-\gamma_l)n/(l-1)}{n/l}=\frac{1-\gamma_l}{1-(1/l)}$$
from the total number of blocks.  Hence 
we find that the proportion of left out blocks
in the optimal alignment is typically close or below the
following bound:
\begin{equation}
\label{firstinequalityq}
\frac{1-\gamma_l}{1-(1/l)}.
\end{equation}
Let us next find a simple lower bound for $\gamma_l$ which we
can use in expression \ref{firstinequalityq}.
Assume we choose an alignment which leaves out no blocks. The typical
score of such an alignment gives a lower bound for
$\gamma_l$. In this case the common subsequence
defined by such an alignment has its $i$-th block
having length:
$$B_i:=\min\{B_{Xi},B_{Yi}\}.$$
where $B_{Xi}$ (resp. $B_{Yi}$) is the length of the $i$-th block of $X$ (resp. $Y$). Recall that $B_{Xi}$ (resp. $B_{Yi}$) has uniform distribution on the set $\{l-1,l,l
+1\}$. The distribution of the minimum above is as follows:
$$\mathrm{P}(B_i=l-1)=5/9,\mathrm{P}(B_i=l)=3/9,\mathrm{P}(B_i=l+1)=1/9.$$
The expected length is thus:
\begin{equation}\label{expBi}
\mathrm{E}[B_i]=\frac59(l-1)+\frac39l+\frac19(l+1)=l-\frac49.
\end{equation}
Since there are about $n/l$ blocks, the score aligning all the blocks
gives thus about a score of:
$$\frac{n}{l}\cdot \mathrm{E}[B_i]=n\left(1-\frac{4}{9l}\right),$$
so that we obtain:
$$\gamma_l\geq \left(1-\frac{4}{9l}\right).$$
The last inequality together with
the bound \ref{firstinequalityq} implies that the proportion
of left out blocks should typically not be much above the following bound:
 \begin{equation} \label{secondinequalityq}
\frac{1-(1-(4/9l))}{1-(1/l)}=\frac{4/9}{l-1}
\end{equation}

\noindent Another similar approach is to get a lower bound for $\gamma_l$ by simulations.
As a matter of fact we have for any $n$ that $\mathrm{E}[L_n]/n$ is a lower bound
for $\gamma_l$. By Montecarlo we can find an estimate of $\mathrm{E}[L_n]/n$
and a very likely lower bound $\gamma_{lb}$. We then replace
in inequality \ref{firstinequalityq} $\gamma_l$
by $\gamma_{lb}$.
\section{High probability events}\label{details}
Let $\delta>0$ be a parameter not depending on $n$. We will define a number of events related with the combinatorial properties of the optimal alignments, called $C^n$, $D^n(\delta)$, $G^n(\delta)$ and $J^n(\delta)$. In the following we will prove that these events have high probability for $n$ large. By high probability, we mean a quantity which is negatively exponential close to one in $n$. It will turn out that this is true for the above events for any parameters $\delta>0$ not depending on
$n$. Also we will prove that $F^n(q)$ has high probability for $n$ large in the same sense as above but restricted to some values of $q$. All the missing proofs (omitted for shortness reasons) can be found with details in \cite{FelipePhD}.

\vspace{12pt}
\noindent A very useful tool we often use is the Azuma-Hoeffding theorem. The following is a version of it for martingales (for a proof see \cite{GS}):
\begin{theorem}(Hoeffding's inequality)\label{hoeffding}
Let $(V,\mathfrak{F})$ be a martingale, and suppose that there exists a sequence $\mathfrak{a}_1,\mathfrak{a}_2,\cdots$ of real numbers such that
\[ \mathrm{P}(|V_n-V_{n-1}| \le \mathfrak{a}_n)=1\]
for all $n$. Then:
\begin{equation}
\mathrm{P}(|V_n-V_0| \ge v) \le 2 \exp\Big\{ -\frac{1}{2} v^2 \Big/ \sum_{i=1}^n \mathfrak{a}^2_i \Big\}\label{hoeffdingIneq}
\end{equation}
for every $v>0$.
\end{theorem}
We also will use a corollary of the above theorem, for some intermediate bounds:
\begin{corollary}\label{hoeffdingco}
Let $a>0$ be constant and $V_1,V_2,\dots$ be an i.i.d sequence of random bounded variables such that:
\[ {\rm P}(|V_i-{\rm E}[V_i]| \le a)=1\]
for every $i=1,2,\dots$ Then for every $\Delta>0$, we have that:
\begin{equation}\label{corAH}
{\rm P}\left(\,\left|\frac{V_1+\dots+V_n}{n}-{\rm E}[V_1] \right|\ge \Delta\right) \le 2\exp\left( -\frac{\Delta^2}{2a^2}\cdot n\right)
\end{equation}
\end{corollary}
\subsection{Number of blocks as renewal process}\label{nblocksRP}
For $k>0$ let us define the sum of the length of the first $k$ blocks in $X$ as:
\[ S^X_k = B_{X1} + \dots + B_{Xk}\]
Let us define the number of blocks used in a sequence of length $t$ in $X$ as:
\begin{equation}
N^X_t = \max \{ k>0 : S^X_k \le t\} \label{defN}
\end{equation}
Note that there might be at the end of $X$ a block which has length smaller than $l-1$. Since this is at most one block it plays little role and we will not mention it every time, only when it is relevant (the same will apply to $Y$ in what follows).
\\
\\
\noindent Due to the standard theory of renewal processes, for every $k,t>0$ the following relation holds between the two random variables defined above:
\begin{eqnarray}\label{fundrt}
N^X_t \ge k \Leftrightarrow S_k^X \le t.
\end{eqnarray}
In the same way we define for $Y$ the same variables as before:
\begin{eqnarray}
S^Y_k &=& B_{Y1} + \dots + B_{Yk} \nonumber \\
N^Y_t &=& \max \{ k>0 : S^Y_k \le t\} \nonumber
\end{eqnarray}
where still the relation $N^X_t \ge k \Leftrightarrow S_k^X \le t$, for every $k,t>0$ holds true. 
\\
\\
\noindent Let $C^n$ be the event that the number of blocks in $X$ and in $Y$ lies in the interval 
\[I_n:=\left[\frac{n}{l}-n^{0.6},\frac{n}{l}+n^{0.6}\right].\]
\begin{lemma} There exists a constant $b_1>0$ depending on $l$ such that:
\[ \mathrm{P}(C^{nc}) \le 8\,e^{-b_1\cdot n^{0.2}}\]
for every $n>0$ large enough.
\end{lemma}
\vspace{12pt}{\bf Proof.} 
It is easy to see that:
\begin{equation}
C_n = \{ N^X_n \in I_n\} \cap \{ N^Y_n \in I_n\} \label{sumaI}
\end{equation}
It is sufficient to compute directly $\mathrm{P}(\{N^X_n \in I_n\}^c)$:
\begin{equation}\label{e1}
\mathrm{P}(\{N^X_n \in I_n\}^c) \le \mathrm{P}\left(N^X_n \le \frac{n}{l} - n^{0.6}\right) + \mathrm{P}\left(N^X_n \ge \frac{n}{l}+n^{0.6}\right)
\end{equation}
Now let us compute each expression separately. Let $m_1:=\left\lceil\frac{n}{l}-n^{0.6}\right\rceil$ be an auxiliar variable. We have at the beginning:
\begin{eqnarray}
\mathrm{P}\left(N^X_n \le \frac{n}{l} - n^{0.6}\right) &\le&\mathrm{P}(N^X_n \le m_1) \nonumber \\
\mbox{(by using $N^X_t \ge k \Leftrightarrow S^X_k \le t$)}&=&{\rm P}\left(S^X_{m_1}\ge n\right)\nonumber\\
&=&{\rm P}\left(\frac{S^X_{m_1}}{m_1}-l \ge \frac{n}{m_1}-l\right) \nonumber\\
\mbox{(by \ref{corAH} with $\mathrm{P}(|B_{X1} - l| \le 1)=1$)}&\le&2\exp\left( -\frac{m_1}{2}\left( \frac{n}{m_1}-l\right)^2\right)\label{a1}
\end{eqnarray}
Now we need to bound $m_1$ in order to get the right order for moderate deviations. Let us start looking at the following:
\begin{eqnarray}
\left(\frac{n}{m_1}-l\right)^2&\ge&l^2\left(\frac{n}{n-ln^{0.6}+l}-1\right)^2\mbox{, by using $m_1 \le \frac{n}{l}-n^{0.6}+1$}\nonumber\\
&\ge&l^2\left( \frac{1}{1-\frac{l}{n^{0.4}}+\frac{l}{n}}-1\right)^2\nonumber\\
&\ge&l^4\left( \frac{1}{n^{0.4}}-\frac{1}{n}\right)^2\left( \frac{1}{1-\frac{l}{n^{0.4}}+\frac{l}{n}}\right)^2\nonumber\\
&\ge&\frac{l^4}{n^{0.8}}\left( 1-\frac{1}{n^{0.6}}\right)^2\left( \frac{1}{1-\frac{l}{n^{0.4}}+\frac{l}{n}}\right)^2\label{a2}
\end{eqnarray}
We have:
\[ \lim_{n \to \infty}\left( 1-\frac{1}{n^{0.6}}\right)^2\left( \frac{1}{1-\frac{l}{n^{0.4}}+\frac{l}{n}}\right)^2 =1 > \frac{1}{4} \]
Hence for $n$ large enough, the expression on the right hand side of \ref{a2} is larger than $l^4/(4n^{0.8})$ so that:
\begin{equation}\label{a2n}
\left(\frac{n}{m_1}-l\right)^2 \ge \frac{l^4}{4n^{0.8}}
\end{equation}
Also, for $n>0$ large enough we can take:
\begin{equation}
m_1=\left\lceil \frac{n}{l}-n^{0.6}\right\rceil \ge\frac{n}{l}-n^{0.6}+1 \ge \frac{n}{2l}+1 =\frac{n}{2l}\left( 1+\frac{2l}{n}\right) \ge \frac{n}{4l}\label{a3}
\end{equation}
Finally we can use \ref{a2n}, \ref{a3} in \ref{a1} to get:
\begin{eqnarray}
\mathrm{P}\left(N^X_n \le \frac{n}{l} - n^{0.6}\right) &\le&2\exp\left( -\frac{m_1}{2}\left( \frac{n}{m_1}-l\right)^2\right) \nonumber\\
&\le&2\exp\left( -\frac{m_1}{2}\cdot\frac{l^4}{4n^{0.8}}\right) \nonumber\\
&\le&2\exp\left( -\frac{n}{4l}\cdot\frac{l^4}{8n^{0.8}}\right) \nonumber\\
&\le&2\exp\left( -\frac{l^3}{32}\cdot n^{0.2}\right) \label{el1}
\end{eqnarray}
for $n>0$ large enough. For the other term, let $m_2:=\left\lfloor\frac{n}{l}+n^{0.6}\right\rfloor$ be an auxiliar variable and do the same as before. We have at the begining:

\begin{eqnarray}
\mathrm{P}\left(N^X_n \ge \frac{n}{l} + n^{0.6}\right) &\le&\mathrm{P}(N^X_n \ge m_2) \nonumber \\
\mbox{(by using $N^X_t \ge k \Leftrightarrow S^X_k \le t$)}&=&{\rm P}\left(S^X_{m_2}\le n\right)\nonumber\\
&=&{\rm P}\left(\frac{S^X_{m_2}}{m_2}-l \le \frac{n}{m_2}-l\right) \nonumber\\
\mbox{(by \ref{corAH} with $\mathrm{P}(|B_{X1} - l| \le 1)=1$)}&\!\!\!\!\!\!\le&\!\!\!\!\!\!2\exp\left( -\frac{m_2}{2}\left( \frac{n}{m_2}-l\right)^2\right)\label{aa1}
\end{eqnarray}
Now we need to bound $m_2$ in order to get the right order for moderate deviations. Let us start looking at the following:
\begin{eqnarray}
\left(\frac{n}{m_2}-l\right)^2&\ge&l^2\left(\frac{n}{n+ln^{0.6}}-1\right)^2\mbox{, by using $m_2 \le \frac{n}{l}+n^{0.6}$}\nonumber\\
&\ge&l^2\left( \frac{1}{1+\frac{l}{n^{0.4}}}-1\right)^2\nonumber\\
&\ge&\frac{l^4}{n^{0.8}}\left( \frac{1}{1+\frac{l}{n^{0.4}}}\right)^2\label{aa2}
\end{eqnarray}
where the very last inequality was obtained by assuming $n$ large enough and noticing that:
\[ \lim_{n \to \infty}\left( \frac{1}{1+\frac{l}{n^{0.4}}}\right)^2 \ge \frac{1}{4} \]
Also, for $n>0$ large enough we can take:
\begin{equation}
m_2=\left\lfloor \frac{n}{l}+n^{0.6}\right\rfloor \ge \frac{n}{2l}\label{aa3}
\end{equation}
Finally we can use \ref{aa2}, \ref{aa3} in \ref{aa1} to get:
\begin{eqnarray}
\mathrm{P}\left(N^X_n \ge \frac{n}{l} + n^{0.6}\right) &\le&2\exp\left( -\frac{m_2}{2}\left( \frac{n}{m_2}-l\right)^2\right) \nonumber\\
&\le&2\exp\left( -\frac{m_2}{2}\cdot\frac{l^4}{4n^{0.8}}\right) \nonumber\\
&\le&2\exp\left( -\frac{n}{4l}\cdot\frac{l^4}{4n^{0.8}}\right) \nonumber\\
&\le&2\exp\left( -\frac{l^3}{16}\cdot n^{0.2}\right) \label{el2}
\end{eqnarray}
\\
Then combining \ref{e1}, \ref{el1} and \ref{el2} we obtain:
\[\mathrm{P}(\{N^X_n \in I_n\}^c) \le 4\,\exp\left( -n^{0.2}\cdot\frac{l^3}{32}\right)\]
and by symmetry we finally get:
\[ {\rm P}(C_n^c) \le 8\,\exp\left( -n^{0.2}\cdot\frac{l^3}{32}\right)\]
for every $n>0$ large enough. Taking $b_1=\frac{l^3}{32}>0$ the proof is finished. $\quad\eop$
\subsection{Left out blocks in an optimal alignment}\label{nLOBOA}
Let $F^n(q)$ denote the already defined event that any optimal
alignment of $X$ and $Y$ leaves out at most a proportion $q$
of blocks in $X$ as well as in $Y$. 
\begin{lemma}
For any $q$ satisfying  $q> \frac{4}{9(l-1)}$, we have that there exists $\beta>0$ such that:
$$\mathrm{P}(F^{nc}(q))\leq e^{-\beta n}$$
for all $n$. Note that here $q$ does not depend on $n$ and also $\beta>0$ does not depend on $n$ but on $l$ and $q$.
\end{lemma}
\begin{lemma}
For every $0<\delta<\frac{1}{l}$ there exists a constant $b_3>0$ depending on $\delta$ and $l$ but not on $n$ such that:
\[ \mathrm{P}(K^{nc}(\delta)) \le e^{-b_3 \cdot  n}\]
for $n$ large enough.
\end{lemma}
\subsection{Proportion of blocks in $X$ and $Y$}\label{ppXY}
Let $X^m$ be the sequence $X^\infty$ taken up to the $m$-th block.
Similarly, let $Y^m$ be the sequence $Y^\infty$ taken up to the
$m$-th block. Let $D_m(\delta)$ be the event that the proportions of blocks
in $X^m$ and $Y^m$ of length $l-1$, $l$ and $l+1$ are not further 
from $1/3$ than $\delta$. Let $D^n(\delta)$ be the event:
\[D^n(\delta)=\bigcap_{m\in I_n}D_m(\delta)\]
where we defined the interval $\,I_n=\left[\,n/l-n^{0.6}\,,\,n/l+n^{0.6}\,\right]$.
\begin{lemma} For every $\delta>0$ we have that:
\[ \mathrm{P}(D^{nc}(\delta)) \le 2n^{0.6}\left(\frac{1}{1+3\delta}\right)^{n\frac{(1+3\delta)}{2l}}\]
for $n$ large enough.
\end{lemma}
\subsection{Number of blocks for an optimal alignment}\label{nboa}
Recall that $N_n^X$ (resp. $N_n^Y$) is the number of blocks in $X$ (resp. in $Y$) having lenghts in $\{l-1,l,l+1\}$ as in expression \ref{defN}. Let $G^n(\delta)$ be the event that the following inequality holds:
$$\frac{N^Y_n}{N_n^X}\leq 1+\delta.$$
\begin{lemma} For every $\delta >0$ there exist two constantsl $b_4,b_5>0$ depending on $l$ and on $\delta$ such that:
\[\mathrm{P}(G^{nc}(\delta)) \le 2e^{-b_4 \cdot n^{0.2}}+2e^{-b_5\cdot n^{0.2}}\]
for every $n$ large enough.
\end{lemma}
\subsection{Cut blocks at the end}\label{cutbEnd}
Let $J^n(\delta)$ denote the event that the proportion of left out blocks at the end of $X$ or $Y$ in any optimal alignment is at most a proportion $\delta$ of the total number of blocks in each of these sequences. As all events before, we want to prove that $J^n(\delta)$ has high probability to happen for every $\delta>0$ provided $n$ is large enough. We need an extra definition and a previous lemma in order to show the high probability of $J^n(\delta)$. 
\\
\\
\noindent For an integer number $s \in [1,n]$ we denote:
\begin{equation}
L_1^s := |{\rm LCS}(X_1X_2 \cdots X_s,Y_1Y_2\cdots Y_n)|\label{lcss}
\end{equation}
\begin{lemma}\label{amv}
Given $\delta>0$, there exists a constant $c_*>0$ not depending on $n$ but on $\delta$ such that:
\begin{equation}\label{ecavm}
{\rm E}[L_n-L_1^{n-\delta n}] \ge c_*\cdot n
\end{equation}
for every $n>0$ large enough.
\end{lemma}
\vspace{12pt}{\bf Proof.} 
Given $n>0$ and $t \in [-1,1]$ let us define the number $\gamma(t,n)>0$ as follows:
\[\gamma(t,n):=\frac{{\rm E}[\,|{\rm LCS}(X_1\cdots X_{n+nt},Y_1\cdots Y_{n-nt})|\,]}{n}\]
This number $\gamma(t,n)$ is a kind of extension for the Chvatal-Sankoff constant $\gamma$ (see \cite{CS}), or more precisely in the case of our paper an extension of $\gamma_l$ defined as in expression \ref{gammal}. An extended motivation for this definition can be found in \cite{AMV}. For any fixed $t \in[-1,1]$ it is known that $\gamma(t,n)$ converges as $n\to \infty$ (see \cite{Alexander} or \cite{AMV}), let us denote that limit by
\[ \gamma(t):=\lim_{n\to \infty} \gamma(t,n).\]
The speed of convergence to that limit is also known due to theorem 2.1 in \cite{Alexander}. This theorem says that there exists $\theta_1>0$ a constant not depending on $n$ such that:
\begin{equation}\label{e1amv}
|\gamma(t,n)-\gamma(t)| \le \frac{\theta_1 \ln(n)}{\sqrt{n}}
\end{equation}
for any fixed $t \in[-1,1]$ provided $n>0$ is large enough.
On the other hand, it is known that the map $t \in[-1,1] \mapsto \gamma(t)\in[0,1]$ is concave and symmetric in the origin (see \cite{AMV}). Hence, for every $t \in [-1,1]$ we have
\begin{equation}\label{e2amv}
\gamma(0) \ge \gamma(t)
\end{equation}
Let us set an auxiliar variable $n^*$ as follows:
\[n^*:= n\left(1-\frac{\delta}{2}\right)\]
Note that with the last definition, the inequality
\begin{equation}\label{especialamv}
\frac{n\ln(n)}{\sqrt{n}}=\sqrt{n}\,\ln(n)\ge\sqrt{n^*}\,\ln(n^*)= \frac{n^*\ln(n^*)}{\sqrt{n^*}}
\end{equation}
holds due to $\sqrt{\cdot}$ and $\ln(\cdot)$ being increasing functions. By using the previous definitions, the inequality for the speed of convergence \ref{e1amv}, the concave inequality \ref{e2amv} and inequality \ref{especialamv} (following this order), we can write:
\begin{eqnarray}
{\rm E}[L_n-L_1^{n-\delta n}] &=& n\,\gamma(0,n)-n^*\gamma\left(t^*,n^*\right) \nonumber \\
&\ge& n\left( \gamma(0) - \frac{\theta_1\ln(n)}{\sqrt{n}}\right)-n^*\left( \gamma(t^*)+\frac{\theta_1\ln(n^*)}{\sqrt{n^*}}\right) \nonumber\\
&\ge& n\left( \gamma(0) - \frac{\theta_1\ln(n)}{\sqrt{n}}\right)-n^*\left( \gamma(0)+\frac{\theta_1\ln(n^*)}{\sqrt{n^*}} \right) \nonumber\\
&=&(n-n^*)\gamma(0)-\theta_1\left( \frac{n\ln(n)}{\sqrt{n}}+\frac{n^*\ln(n^*)}{\sqrt{n^*}}\right) \nonumber \\
&\ge&(n-n^*)\gamma(0)-2\theta_1\frac{n\ln(n)}{\sqrt{n}}\nonumber\\
&=&\frac{n\,\delta\,\gamma(0)}{2}-2\theta_1\frac{n\ln(n)}{\sqrt{n}}\nonumber\\
&=&\left(\frac{\delta\,\gamma(0)}{2}-2\theta_1\frac{\ln(n)}{\sqrt{n}}\right)n\nonumber\\
&\ge&\frac{\delta\,\gamma(0)}{4}\,n\label{lastamv}
\end{eqnarray}
where the very last inequality above holds for $n$ large enough, since
\[ \lim_{n\to\infty}\left(2\theta_1\frac{\ln(n)}{\sqrt{n}}\right) =0 < \frac{\delta\,\gamma(0)}{4}\]
To finish the proof we take $c_*=\frac{\delta\,\gamma(0)}{4}$.$\quad\eop$

\vspace{12pt}
\noindent Now comes the main result of this section which establishes the high probability of the event $J^n(\delta)$:
\begin{proposition}\label{remarkInc}
For every $\delta>0$, there exists a constant $\theta>0$ not depending on $n$ but on $\delta$ such that:
\[ {\rm P}(J^{nc}(\delta)) \le 2e^{-\theta\cdot n}\]
for every $n>0$ large enough.
\end{proposition}
\vspace{12pt}{\bf Proof.} 
With the notation as in \ref{lcss} we write:
\begin{eqnarray}
{\rm P}(J^{nc}(\delta))&\le& 2\,{\rm P}(\,|{\rm LCS}(X_1\cdots X_{n-\delta n},Y_1 \cdots Y_n)|-L_n \ge 0\,)\nonumber\\
&=&2\,{\rm P}(\,L_1^{n-\delta n}-L_n \ge 0\,)\nonumber \\
&=&2\,{\rm P}(\,L_1^{n-\delta n}-L_n -{\rm E}[\, L_1^{n-\delta n}-L_n\,]\ge {\rm E}[\, L_n-L_1^{n-\delta n}\,]\,)\label{v1}
\end{eqnarray}
Let us define
\[M_n(\delta):=L_1^{n-\delta n}-L_n -{\rm E}[\, L_1^{n-\delta n}-L_n\,]\]
It is not difficult to see that $M_n(\delta)$ is a martingale with respect to the filtration $\mathfrak{F}_n=\sigma\{(X_k,Y_k):k \le n \}$ and that $M_0=0$. The following inequality also holds:
\[ |M_n(\delta)-M_{n-1}(\delta)| \le 4\]
for $\delta>0$ with probability 1. So, we can use the theorem \ref{hoeffding} (Azuma-Hoeffding inequality for martingales) with $\mathfrak{a}_i=4$ and $v={\rm E}[\, L_n-L_1^{n-\delta n}\,]$ to estimate:
\begin{eqnarray}
{\rm P}(\,L_1^{n-\delta n}-L_n -{\rm E}[\, L_1^{n-\delta n}-L_n\,]\ge v\,) &\le &2\exp\left(-\frac{v^2}{2\cdot 4n}\right)\nonumber\\
&=&2\exp\left(-\frac{n}{8}\left(\frac{{\rm E}[\, L_n-L_1^{n-\delta n}\,]}{n}\right)^2\right) \nonumber \\
\mbox{(by \ref{ecavm} and $c_*$ from lemma \ref{amv})}&\le& 2\exp\left(-\frac{c_*^2}{8}\cdot n\right) \nonumber
\end{eqnarray}
Taking $\theta=\frac{c_*^2}{8}>0$ finishes the proof.
$\quad\eop$
\subsection{Optimal events}\label{optimalEvents}
In theorem \ref{minimizingtheo}, $q$ represents the proportion of left out blocks in $X$ and in $Y$. In reality, typically, the proportion of left out blocks in $X$ will not be exactly equal to the proportion of left out blocks in $Y$. Because of this, $q_1$ will designate the proportion of left out blocks in $X$ and $q_2$ will designate the proportion of left out blocks in $Y$. We will have that $q_1$ can be made as close to $q_2$ as we want to by taking a large $n$. Now we need to rewrite all our conditions as in theorem \ref{minimizingtheo} with $q_1$ and $q_2$ instead of $q$.
\\
\\
\noindent Let us define the following events:
\begin{itemize}
\item{} Given any $m_1,m_2,q_1,q_2$, let $E_{m_1,m_2,q_1,q_2}(\epsilon)$ denote the event that there
is no optimal alignment of $X^{m_1}$ with $Y^{m_2}$ leaving out
a proportion of $q_1$ blocks in $X^{m_1}$ and a proportion
of $q_2$ blocks in $Y^{m_2}$ and such that:
\begin{equation}
\label{fundamentalepsilon}
H(q_1)+H(q_2)+(1-\max\{q_1+3q_2,3q_1+q_2\})
\left(\ln(1/9)+H(p)\right)\leq -\epsilon.
\end{equation}

\item{} Let $E^n(\epsilon)$ be the event :
\begin{equation}\label{defofEn}
E^n(\epsilon)=\bigcap_{m_1,m_2\in I^n,q_1,q_2}E_{m_1,m_2,q_1,q_2}(\epsilon).
\end{equation}
\end{itemize}

\noindent If $\delta$ designates the difference between $q_1$ and $q_2$, then note that the system
\begin{eqnarray}
\min\left[\,\,\frac{p_{l-1,l}+p_{l-1,l+1}}{p_{l-1,l-1}+p_{l-1,l}+p_{l-1,l+1}}\left(1-\frac{3q_1}{(1/3)-\delta}  \right)+\right.&&\nonumber\\
\left.-\left(1-\frac {\delta+2(q_1-\delta)}{(1/3+\delta)}\right)\frac{p_{l+1,l+1}}{p_{l+1,l-1}+p_{l+1,l}+p_{l+1,l+1} }-q_2\cdot\frac{1+\delta}{(1/3)-\delta}\,\,\right]\nonumber
\end{eqnarray}
\vspace{-12pt}
\begin{eqnarray}
|q_1-q_2| &\leq&2\delta \nonumber \\
H(q_1)+H(q_2)+(1-\max\{q_1+3q_2,3q_1+q_2\})\left(\ln(1/9)+H(p)\right)&\geq& 0 \nonumber\\ \label{christ}
\end{eqnarray}
converges to the conditions in theorem \ref{minimizingtheo}
when $\delta$ goes to zero (when $q_1$ is as close to $q_2$ as we want to by taking a large $q_1$). Note also that replacing $q_1$ and $q_2$ by $q$ and taking $\delta=0$ in the minimized function and in the last inequality of \ref{christ}, they become equal to \ref{expressionI} respectively \ref{entropycondition}. If the minimizing problem in theorem \ref{minimizingtheo} has a strictly positive solution $2\epsilon$ and if expression \ref{expressionI} is less than $\epsilon_1$, this implies that \ref{entropycondition} must be smaller than a $-\epsilon_2$ for $\epsilon_2>0$ (we are assuming that \ref{conditionq}, \ref{conditionpij1} and \ref{entropycondition} hold). The next lemma shows that the same holds true for the system \ref{christ} if we take $\delta$ small enough.

\begin{lemma}\label{approachOpt}
Assume there exists $0<q_0<(1/3)$ and
 $\epsilon_1>0$ such that for all $\{p_{ij}\}_{i,j}$ and $q\in[0,q_0]$
 satisfying all the conditions \ref{conditionq}, 
\ref{conditionpij1} and \ref{entropycondition} in theorem \ref{minimizingtheo},
we have that
expression \ref{expressionI} 
is larger or equal to $2\epsilon_1$ (in other words, the condition that the minimizing problem in theorem \ref{minimizingtheo} has a strictly positive solution $2\epsilon_1$ is satisfied). Then, we have that
there exists $\epsilon_2>0$ and $\delta_0>0$ such that
for all $\{p_{ij}\}_{i,j\in \{l-1,l,l+1\}}$ and $q_1,q_2\in[0,q_0]$ and 
$\delta\in[0,\delta_0]$ satisfying
\ref{conditionq} and \ref{conditionpij1},
we have that if $|q_1-q_2|\leq 2\delta_0$ and if
\begin{eqnarray}
\frac{p_{l-1,l}+p_{l-1,l+1}}{p_{l-1,l-1}+p_{l-1,l}+p_{l-1,l+1}}\left(1-\frac{3q_1}{(1/3)-\delta_0}  \right)&&+\nonumber\\
-\left(1-\frac {\delta_0+2(q_1-\delta_0)}{(1/3+\delta_0)}\right)\frac{p_{l+1,l+1}}{p_{l+1,l-1}+p_{l+1,l}+p_{l+1,l+1} }-q_2\cdot\frac{1+\delta_0}{(1/3)-\delta_0}&\leq& \epsilon_1\nonumber \\\label{expressionI*} 
\end{eqnarray}
then
\begin{equation}
\label{Hq1+Hq2}
H(q_1)+H(q_2)+(1-\max\{q_1+3q_2,3q_1+q_2\})
\left(\ln(1/9)+H(p)\right)\leq -\epsilon_2
\end{equation}
\end{lemma}
{\bf Proof.}  We are going to do the proof by {\it reductio ad absurdum} (reduction to the absurd). Assume for this that
for all $(p_{ij})_{i,j\in I}$ and $q\in[0,q_0]$
 satisfying all the conditions \ref{conditionq}, 
\ref{conditionpij1} and \ref{entropycondition} in theorem \ref{minimizingtheo}
we have that
expression \ref{expressionI} 
is larger equal to $2\epsilon_1$. Assume that the rest of the lemma would not
hold. Then for every $\delta>0$ (as small as we want)
we could find a vector $\vec{p}$:
$$\vec{p}:=(p_{l-1,l-1},p_{l-1,l},\ldots,p_{l+1,l+1},q_1,q_2,\delta)$$
such that the components satisfy $|q_1-q_2|\leq \delta$,
and the components of $\vec{p}$ satisfy \ref{conditionq}, 
\ref{conditionpij1} whilst inequality \ref{expressionI*} is satisfied
and we can take the expression
\begin{equation}\label{entropyconditionq1q2}
H(q_1)+H(q_2)+(1-\max\{q_1+3q_2,3q_1+q_2\})\left(\ln(1/9)+H(p)\right)
\end{equation}
as close to zero as we want.  Hence there exists
a sequence $\vec{p}_1,\vec{p}_2,\ldots,\vec{p}_t,\ldots$ of such vectors
with notation:
$$\vec{p}(t):=(p_{l-1,l-1}(t),p_{l-1,l}(t),\ldots,p_{l+1,l+1}(t),q_1(t),q_2(t),\delta(t))$$
so that for each $t\in\mathbb{N}$ the vector $\vec{p}(t)$ satisfies all the conditions \ref{conditionq}, 
\ref{conditionpij1} and \ref{expressionI*},
whilst 
$$\lim_{t\rightarrow\infty}|q_1(t)-q_2(t)|=0$$
and expression \ref{entropyconditionq1q2} converges to zero as
$t$ goes to infinity.\\
The vectors $\vec{p}(t)$ are contained in a bounded domain
and hence in a compact domain. This implies that there exists
a converging subsequence. Hence there exists an increasing map
$\pi:\mathbb{N}\rightarrow\mathbb{N}$ so that
$\vec{p}(\pi(t))$ converges as $t$ goes to infinity.
Let the limit be denoted by
$$\vec{p}:=(p_{l-1,l-1},p_{l-1,l},
\ldots,p_{l+1,l+1},q_1,q_2,0).$$
We have that $q_1=q_2$, so let us denote $q_1=q_2=q$.
We find that our limit satisfies all the conditions \ref{conditionq},  \ref{conditionpij1}. 
Furthermore, at the limit expression \ref{entropyconditionq1q2} becomes equal to zero.
Replacing then $q_1$ and $q_2$ in \ref{entropyconditionq1q2} by $q$,
we find that condition \ref{entropycondition} is satisfied.
finally since for our sequence $\vec{p}(\pi(t))$, we have
that \ref{expressionI*} is satisfied, by continuity it must
also be satisfied for the limit. Hence, noting that at 
the limit $q_1=q_2=q$ and $\delta=0$, we get
that expression \ref{expressionI} is less or equal to $\epsilon_1$.
This contradicts our assumption, since our limit 
vector satisfies all the conditions\ref{conditionq}, 
\ref{conditionpij1} and \ref{entropycondition}  and should thus
have expression 
\ref{expressionI} 
 larger equal to $2\epsilon_1$. Hence, we have that when all the conditions
\ref{conditionq} and \ref{conditionpij1}, 
are satisfied and when $\delta$ goes to zero then
expression \ref{expressionI*} should be bounded away from zero.
This means that for $\delta>0$ small enough and when all the conditions
\ref{conditionq} and \ref{conditionpij1} are
satisfied, we have that there exists $\epsilon_2>0$ so that
\ref{expressionI*} is less or equal to $-\epsilon_2$.
$\quad\eop$

\vspace{12pt}
\noindent Let us now show that event $E_{m_1,m_2,q_1,q_2}(\epsilon)$ holds with high probability.

\begin{lemma}\label{boundEepsilon}
Assume that there exists $0<q_0<(1/3)$ and $\epsilon_1>0$ such that for all $\{p_{ij}\}_{i,j}$ and $q\in[0,q_0]$ satisfying all the conditions \ref{conditionq}, \ref{conditionpij1} and \ref{entropycondition} in theorem \ref{minimizingtheo}, we have that expression \ref{expressionI} is larger or equal to $2\epsilon_1$. Then, for every $\epsilon>0$ there exist a polynomial $w(n)>0$ and a constant $\vartheta >0$, both only depending on $l$ such that: 
\[\mathrm{P}(E^{nc}(\epsilon)) \leq w(n) e^{-\vartheta\cdot  n}\]
for every $n$ large enough.
\end{lemma}
{\bf Proof.} 
Let $\vec{a}$ denote an alignment of the $X^{m_1}$ and $Y^{m_2}$.
Hence $a$ consists of two binary vectors 
$\vec{a}=(\vec{a}_X,\vec{a}_Y)$ the first one having 
length $m_1$ and the second one having length $m_2$.

\noindent Hence $\vec{a}_x\in \{0,1\}^{m_1}, \vec{a}_Y\in\{0,1\}^{m_2}$ when the $i$-th entry $a_{Xi}$ of $\vec{a}_X$ is 
a $1$ that means that the $i$-th block
of $X_{m_1}$ is discarded (entirely aligned with gaps)by 
the alignment $\vec{a}$, otherwise the $i$-th block
of $X^{m_1}$ is not discarded. Similarly when
$a_{Yi}=1$ then the $i$-th block of $Y^{m_1}$ is
discarded.  Here we use the same way of defining
alignment as explained before in the first section:
we specify which blocks get entirely discarded
and then align the rest block by block.
Doing so and assuming that the alignment
$\vec{a}$ is not random, we get that the
aligned block pairs are i.i.d.. For the lengths of
aligned block pairs we have nine possibilities
each having the same probability. Hence, given the alignment
$a$, the empirical frequencies of the aligned block
pair lengths is simply a multinomial distribution.
Let $p=\{p_{ij}\}_{i,j\in\{l-1,l,l+1\}}$ be a (non-random)
probability distribution. Let  $E_{a}(p)$ denote the event that the
empirical distribution of the aligned block pairs
by the alignment $a$ is not $p$.

\noindent From what we said we have that
the probability $\mathrm{P}(E_a^c(p))$ is equal
to the probability that a 9-nomial variable with
parameter $m^*$ and all probability parameters 
equal to $1/9$ gives the frequencies
given by $p$. Here $m^*$ designates the number of aligned
block pairs by $a$. Hence, we get
\begin{equation}\label{methodtypes1}
\mathrm{P}(E_a^c(p))=\binom {m^*}{m^*p_{l-1,l-1}\;\;m^*p_{l-1,l}\;\ldots\;m^*p_{l+1,l+1}}
\left(\frac19\right)^{m^*}
\end{equation}
where
\[ \binom{a}{a_1\dots a_k}=\frac{a!}{a_1! \cdots a_k!}\]
is the multinomial factorial coefficient. Let us define
\begin{eqnarray}
B(p) &:=&\binom {m^*}{p_{l-1,l-1}m^*\;\;p_{l-1,l}m^*\;\ldots\;p_{l+1,l+1}m^*} \nonumber \\
M(p)  &:=& \prod_{p_i \in \{ p_{l-1,l-1},\dots,p_{l+1,l+1}\}} p_i^{p_i} \nonumber \\
H(p) &:=& \sum_{p_i \in \{ p_{l-1,l-1},\dots,p_{l+1,l+1}\}} p_i \ln (1/p_i) \;=\; \ln\left(\frac{1}{M(p)}\right) \label{entropiaHp}
\end{eqnarray}
note that $B(p)\cdot (M(p))^{m^*}$ is the probability distribution of a multinomial random variable with parameters $m^*$ and vector $(m^*p_{l-1,l-1},\dots,m^*p_{l+1,l+1})$. Hence
\begin{equation}\label{multp}
B(p)\cdot (M(p))^{m^*} \le 1.
\end{equation}
Then, by using \ref{multp} we can bound expression \ref{methodtypes1} as follows:
\begin{eqnarray}
B(p) \left( \frac{1}{9}\right)^{m^*} &=& B(p)\cdot (M(p))^{m^*}\left( \frac{1/9}{M(p)}\right)^{m^*} \nonumber \\
&\le&\left( \frac{1/9}{M(p)}\right)^{m^*} \nonumber \\
&=& \exp\left(\left[\ln\left(\frac{1}{9}\right)+H(p)\right]m^*\right) \label{bound}
\end{eqnarray}
On the other hand, we have at least $(1-\max\{q_1+3q_2,3q_1+q_2\})\min\{m_1,m_2\}$
aligned block pairs. Let us give an intuition for this. There are three situations for aligning a fixed block in $X$ with blocks in $Y$. First, when we align one block in $X$ with one block in $Y$ one to one, the resulting length contributing to the {\rm LCS} is the minimun between their lenghts, so at most if all the blocks of $X$ and $Y$ are aligned one to one then we will have at most a contribution of $\min\{m_1,m_2\}$ aligned blocks pairs. Second, when we align one block in $X$ with several blocks in $Y$ then we at least leave $q_1 \cdot m_1$ blocks in $X$. Third, when know that we cannot align two adjacent blocks in $X$ with the same block in $Y$, then we leave at least $2q_1\cdot m_1$ blocks in $X$ also. In total, in the worse case, looking first at blocks in $X$, we are leaving $(3q_1+q_2)\min\{m_1,m_2\}$ blocks in both sequences $X$ and $Y$. Similarly, but looking first at $Y$, we can leave $(3q_2+q_1)\min\{m_1,m_2\}$ blocks in both sequences $X$ and $Y$. Finally, at least we have $(1-\max\{q_1+3q_2,3q_1+q_2\})\min\{m_1,m_2\}$ aligned block pairs due to the considerations above.

\noindent Since $m_1,m_2\in I_n$, this gives the lower bound for $m^*$
\begin{equation}
\label{lowerboundm^*}
m^*\geq (1-\max\{q_1+3q_2,3q_1+q_2\})\;\cdot\;
((n/l)-n^{0.6}).
\end{equation}
and hence together with the bound \ref{bound}, we obtain
\begin{equation}
\label{oneinequalitywithEa}
\mathrm{P}((E_a^c(p))\leq
\exp\left(\;(\ln(1/9)+H(p))(1-\max\{q_1+3q_2,3q_1+q_2\})
((n/l)-n^{0.6})\;\right)
\end{equation}
Let $\mathcal{A}_{m_1,m_2,q_1,q_2}$ denote the set of all alignments 
aligning $X^{m_1}$ with $Y^{m_2}$ and leaving out
a proportion of $q_1$ blocks in $X^{m_1}$
and a proportion of $q_2$ blocks in $Y^{m_2}$.
In other words,  the set $\mathcal{A}_{m_1,m_2,q_1,q_2}$
is the set of all elements 
$\vec{a}=(\vec{a}_X,\vec{a}_Y)$
of $\{0,1\}^{m_1}\times\{0,1\}^{m_2}$
for which $|\vec{a}_{X}|=q_1m_1$ 
and $|\vec{a}_Y|=q_2m_2$.\\
Let $\mathcal{P}_{\epsilon,q_1,q_2}$ 
denote the set of those distributions $p$ 
(for aligned block pairs, hence
on the space $\Omega=\{(l-1,l-1),(l-1,l),\ldots,(l+1,l+1)\}$)
for which inequality 
\ref{fundamentalepsilon} is satisfied and which are possible
in our case. Before we continue with the proof, let us look at an example:
\begin{example}
Assume we look at binary strings of length $5$.
Then there can be 0,1,2,3,4 or 5 ones. Hence, the empirical 
distribution for side one when we flip a coin exactly five times
can only be $0$, $20\%$, $40\%$, $60\%$, $80\%$ or $100\%$.
In general for a string of length $n$ and $k$ symbols,
there are no more than $(n+1)^{k-1}$ possible empirical distributions (see \cite{Dembo}, Lemma 2.1.2 (a)). 
In the case above we have an empirical distribution
for $m^*$ aligned block pairs. For each block pairs there
are $9$ possibilities. Hence, there are no more than $(m^*+1)^8$
possible empirical distributions. However $m^*$ is not known.
It could potentially take on any value between $1$ and
$(n/l)+n^{0.6}$. Hence, we find that for the number of 
empirical distributions we need to consider the following upper
bound:
$$((n/l)+n^{0.6}) \cdot ((n/l)+n^{0.6}+1)^8 \le ((n/l)+n^{0.6}+1)^9.$$
\end{example}

\vspace{12pt}
\noindent Let us continue with the proof. We have that: 
$$\bigcap_{a\in\mathcal{A}_{m_1,m_2,q_1,q_2},\mathcal{P}_{\epsilon,q_1,q_2}}E_a(p)=
E_{m_1,m_2,q_1,q_2}(\epsilon)$$
and hence:
\begin{equation}
\label{Em1m2}
\mathrm{P}(E_{m_1,m_2,q_1,q_2}^c(\epsilon)) \leq \sum_{a\in\mathcal{A}_{m_1,m_2,q_1,q_2}\,,\,p \in \mathcal{P}_{\epsilon,q_1,q_2}} \mathrm{P}(E_a^c(p)).
\end{equation}
By using
\ref{oneinequalitywithEa}, the inequality \ref{Em1m2} above becomes:
{\footnotesize \[\mathrm{P}(E_{m_1,m_2,q_1,q_2}^c(\epsilon)) \le  \!\!\!\!\!\!\!\!\!\!\!\!\!\!\!\!\!\!\sum_{a\in\mathcal{A}_{m_1,m_2,q_1,q_2}\,,\,p \in \mathcal{P}_{\epsilon,q_1,q_2}} \!\!\!\!\!\!\!\!\!\!\!\!\!\!\!\!\!\!\exp\left(\; (\ln(1/9)+H(p))(1-\max\{q_1+3q_2,3q_1+q_2\}) ((n/l)-n^{0.6})\;\right).\]}
Note that the number of alignment considered in the sum on the right hand side of the last inequality above
can be bound as follows:
\begin{eqnarray}
|\mathcal{A}_{m_1,m_2,q_1,q_2}| &=& { m_1\choose q_1 m_1 (1-q_1)m_1}{ m_2\choose q_2 m_2 (1-q_2)m_2} \nonumber \\
& \le & \left( \frac{1}{q_1^{q_1}(1-q_1)^{1-q_1}}\right)^{m_1} \left( \frac{1}{q_2^{q_2}(1-q_2)^{1-q_2}}\right)^{m_2} \nonumber \\
&=&\exp(H(q_1)m_1 + H(q_2)m_2) \nonumber \\
& \le & \exp((H(q_1) + H(q_2))m^*) \nonumber \\
& \le & \exp\left(\;(H(q_1)+H(q_2))((n/l)+n^{0.6})\;\right)
\end{eqnarray}
where for $i=1,2$ we denote
\[ H(q_i):=q_i \ln (1/q_i) + (1-q_i)\ln (1/(1-q_i)).\]
The number of distributions in $\mathcal{P}_{\epsilon,q_1,q_2}$ we need to consider is (as explained above) less or equal to $((n/l)+n^{0.6}+1)^9$. Combining all of this we find
that $\mathrm{P}(E_{m_1,m_2,q_1,q_2}^c(\epsilon))$
is less or equal to:
{\footnotesize \[\exp((H(q_1)+H(q_2))((n/l)+n^{0.6})) \cdot b \cdot \exp\left(\;(\ln(1/9)+H(p))(1-\max\{q_1+3q_2,3q_1+q_2\})((n/l)-n^{0.6})\;\right).\]}

\noindent where $b:=((n/l)+n^{0.6}+1)^9$.
In other words, we found that:
{\footnotesize \begin{equation} \label{importantEm}
\mathrm{P}(E_{m_1,m_2,q_1,q_2}^c(\epsilon)) \leq b\exp\left(\frac{n}{l}(H(q_1)+H(q_2)+(\ln(1/9)+H(p))(1-\max\{q_1+3q_2,3q_1+q_2\})+r)\;\right)
\end{equation}}

\noindent where the rest term $r$ is equal to:
$$r=ln^{-0.4}\left(H(q_1)+H(q_2)- (\ln(1/9)+H(p))(1-\max\{q_1+3q_2,3q_1+q_2\})\right)$$
being bounded as follows:
$$ |r| \leq ln^{-0.4}(|H(q_1)|+|H(q_2)|+ (|\ln(1/9)|+|H(p)|))=ln^{-0.4}(3+|\ln(1/9)|).$$
Note that $r$ is bounded from above by a constant times $n^{-0.4}$
where the constant does not depend on $l,q_1,q_2,p$. Hence for $n$ large enough:
\begin{equation}
\label{r}
r\leq \epsilon/2
\end{equation}
Note also that in the sum
\ref{Em1m2}, we only took distributions 
$p \in \mathcal{P}_{\epsilon,q_1,q_2}$ hence satisfying
inequality
\ref{fundamentalepsilon}. This implies that
in the bound \ref{importantEm}, we can assume that
inequality \ref{fundamentalepsilon}
holds. This then implies
\begin{equation}
\label{importantEm2}
\mathrm{P}(E_{m_1,m_2,q_1,q_2}^c(\epsilon)) \leq b\exp\left(\frac{n}{l}(-\epsilon+r)\;\right)
\end{equation}
Assuming now that \ref{r} holds, we obtain:
\begin{equation}
\label{importantEm3}
\mathrm{P}(E_{m_1,m_2,q_1,q_2}^c(\epsilon)) \leq b\exp\left(\frac{n}{l}(-\epsilon/2)\;\right)
\end{equation}
Note that the bound on the right side of the 
last inequality above is negatively exponentially small in $n$,
since $b$ is an expression which is only polynomial in $n$.
Using the equation \ref{defofEn}, we obtain:
$$\mathrm{P}(E^{nc}(\epsilon))\leq\sum_{m_1,m_2\in I^n,q_1,q_2}\mathrm{P}(E_{m_1,m_2,q_1,q_2}^c(\epsilon)).$$
Applying inequality \ref{importantEm3} to the last inequality above, we obtain:
\begin{equation}
\label{defofEn2}
\mathrm{P}(E^{nc}(\epsilon)) \leq \sum_{m_1,m_2\in I^n,q_1,q_2}b\exp\left(\frac{n}{l}(-\epsilon/2)\;\right).
\end{equation}
Note that when $m_1$ is given, the number of possibilities
for the number of left out blocks in $X^{m_1}$ is at most $m_1$.
Hence, for given $m_1$ we have that $q_1$ can take on at most
$m_1$ values. Similarly for given $m_2$ we have that
$q_2$ can take on at most $m_2$ values. But $m_1$ and $m_2$
are less then $(n/l)+n^{0.6}$. Also, both $m_1$ and $m_2$ are
in $I_n$ hence they can take on at most $2n^{0.6}$ values.
This implies that in the sum \ref{defofEn2}, the number of terms
is bound above by the expression:
$$\left((n/l)+n^{0.6}  \right)^24n^{1.2}$$
This upper bound applied to inequality \ref{defofEn2} yields:
\begin{equation}
\label{defofEn3}
\mathrm{P}(E^{nc}(\epsilon))
\leq b\left((n/l)+n^{0.6}  \right)^24n^{1.2}\exp\left(\frac{n}{l}(-\epsilon/2)\;\right).
\end{equation}
which is the negative exponential upper bound we where looking for.
$\quad\eop$
\subsection{Positive expected change in the score}
Let us recall the events that we have proven to have high probability:
\begin{itemize}
\item $C^n$ is the event that the number of blocks in $X$ and in $Y$ lies in the interval
\[ I_n=\left[ \frac{n}{l}-n^{0.6},\frac{n}{l}+n^{0.6}\right].\]
\item $D^n(\delta)$ is the intersection
\[ D^n(\delta)=\bigcap_{m \in I_n}D_m(\delta),\]
where $D_m(\delta)$ is the event that the proportion of blocks in $X^m$ and in $Y^m$ of length $l-1,l$ and $l+1$ are not further from $1/3$ than $\delta$, where $X^m$ (resp. $Y^m$) denotes the sequence $X^\infty$ taken up to the $m$-th block (resp. the sequence $Y^\infty$ taken up to the $m$-th block).
\item $F^n(q)$ is the event that any optimal alignment of $X$ and $Y$ leaves out at most a proportion $q$ of blocks in $X$ as well as in $Y$.
\item $G^n(\delta)$ is the event that the following inequality holds:
\[ \frac{N^Y_n}{N^X_n} \le 1+\delta\]
where $N^X_n$ (resp. $N^Y_n$) is the number of blocks in $X$ (resp. in $Y$) having length in $\{l-1,l,l+1\}$.
\item $E^n(\epsilon)$ is the intersection
\[ E^n(\epsilon) = \bigcap_{m_1,m_2 \in I_n\,;\, q_1,q_2\in[0,1]} E_{m_1,m_2,q_1,q_2}(\epsilon)\]
where $E_{m_1,m_2,q_1,q_2}(\epsilon)$ is the event that there is no optimal alignment of $X^{m_1}$ with $Y^{m_2}$ leaving out a proportion of $q_1$ blocks in $X^{m_1}$ and a proportion of $q_2$ blocks in $Y^{m_2}$ and such that:
\[ H(q_1)+H(q_2)+(1-\max\{q_1+3q_2,q_2,3q_1\})(\ln(1/9)+H(p)) \le -\epsilon_2\]

where $\epsilon_2>0$ depends on $\epsilon, \delta_0$ and $q_0$ and comes from lemma \ref{approachOpt}, $X^{m_1}$ (resp. $Y^{m_2}$) denotes the sequence $X^\infty$ taken up to the $m_1$-th block (resp. the sequence $Y^\infty$ taken up to the $m_2$-th block) and $H(p)$ denotes the entropy as in \ref{entropiaHp} for an alignment.
\end{itemize}

\noindent We can now formulate our combinatorial lemma based on those events:
\begin{lemma}\label{quasitheorem}
Let us consider the constants $q_0, \epsilon_1, \delta_0$ and $\epsilon_2$ from lemma \ref{approachOpt}. Assume that $C^n$, $D^n(\delta_0)$, $F^n(q_0)$, $G^n(\delta_0)$
 and $E^n(\epsilon_2)$
all hold. Then, we have that
$$\mathrm{E}[\tilde{L}_n-L_n|X,Y]\geq \epsilon_1$$
\end{lemma}
\vspace{12pt}{\bf Proof.} For any $x,y\in \{0,1\}^n$ let
$L(x,y)$ denote the length of the {\rm LCS} of $x$ and $y$. 
Let now $x,y\in\{0,1\}$ be any two
realizations so that if $X=x$ and $Y=y$, then
the events $C^n$, $D^n(\delta_0)$, $F^n(q_0)$ and $E^n(\epsilon_2)$
all hold. Let $a$ be a left most optimal alignment of 
$x$ and $y$.  Let $\tilde{x}$ denote the sequence $x$ on which
we performed our random changes. That is $\tilde{x}$ is obtained
by selecting a block of length $l-1$ at random and changing it to length
$l$ and also selecting a block of length $l+1$ at random and reducing it
to length $l$. Let $x^*$ be the sequence we obtain by applying to $x$ 
 only  the first one of the two random changes. That is $x^*$ is obtained
be increasing the length of a randomly chosen block of $x$
of length $l-1$ to length $l$.  So, we start with $x$.
Then we apply the first change and obtain $x^*$. Then
in $x^*$ we choose a block of length $l+1$ at random,
decrease it by one unit to obtain $\tilde{x}$.\\
For all $i,j\in\{l-1,l,l+1\}$, let $p_{ij}$ denote the proportion
of aligned block pairs with lengths $(i,j)$ in the alignment
$a$ of $x$ and $y$. Let $q_1$, resp. $q_2$ denote the proportion
of blocks not aligned by $a$ in $x$, resp. in $y$.
Let $p_{l-1}^I$ denote the proportion of blocks which get aligned by $a$ one block to one block, among all blocks of $x$ of length $1-1$.
Let $p_{l-1}^{II}$ denote the proportion
among all blocks of $x$ of length $l-1$ which are aligned
with several blocks of $y$.  Finally, let $p_{l-1}^{III}$
 denote the proportion among the blocks
of length $l-1$ in $x$ which are left out or are together with other blocks
of $x$ aligned with the same block of $y$. Note that when we increase by
one unit  a block in this third category, then in general the score does
not get any increase.
On the other hand, assume that  the block of $x$ length $l-1$ chosen randomly
and increased by one unit, 
is aligned one block with one block. Then if this chosen block
is aligned with a block of length $l$ or $l+1$ the score is
going to increase. Let $G_{l-1,I}$  be the event that
the block of length $l-1$ chosen is aligned one block with one block.
From what we said it follows that:
$$\mathrm{P}(L(x^*,y)-L(x,y)=1\;|\;G_{l-1,I})\geq 
\frac{ p_{l-1,l}+p_{l-1,l+1} } { p_{l-1,l-1}+p_{l-1,l}+p_{l-1,l+1} }.$$
Note that by only adding a bit the score cannot decrease, so that
the last inequality above means:
\begin{equation}
\label{gl-1I}
\mathrm{E}[L(x^*,y)-L(x,y)\;|\;G_{l-1,I}]\geq 
\frac{p_{l-1,l}+p_{l-1,l+1}}{p_{l-1,l-1}+p_{l-1,l}+p_{l-1,l+1}}.
\end{equation}   
When the block of length $l-1$ chosen and increased is aligned with several
blocks of $y$ at the same time, then we will always observe and increase
of one unit. This yields:
\begin{equation}
\label{gl-1II}
\mathrm{E}[L(x^*,y)-L(x,y)\;|\;G_{l-1,II}]= 
1,
\end{equation}
where $G_{l-1,II}$ denotes the event that the chosen block of length
$l-1$ is aligned with several blocks of $y$. By law of total probability we 
find thus:
\begin{eqnarray}
\mathrm{E}[L(x^*,y)-L(x,y)]&\geq& \mathrm{P}(G_{l-1,I})\frac{p_{l-1,l}+p_{l-1,l+1}}{p_{l-1,l-1}+p_{l-1,l}+p_{l-1,l+1}}+ \mathrm{P}(G_{l-1,II})\nonumber \\
&\geq& \left( 1-\mathrm{P}(G_{l-1,III}) \right) \frac{p_{l-1,l}+p_{l-1,l+1}}{p_{l-1,l-1}+p_{l-1,l}+p_{l-1,l+1}}\nonumber
\end{eqnarray} 
where $G_{l-1,III}$ denotes the event that the block of length
$l-1$ chosen is left out or aligned to the same block of $y$ at the same time as other blocks
of $x$. The last inequality above yields: 
\begin{equation} \label{heini}
\mathrm{E}[L(x^*,y)-L(x,y)]\geq
(1-p_{l-1}^{III})
\frac{p_{l-1,l}+p_{l-1,l+1}}{p_{l-1,l-1}+p_{l-1,l}+p_{l-1,l+1}}.
\end{equation}
Note that the proportion of left out blocks in $x$ is $q_1$. There
can not be two adjacent blocks of $x$ aligned with the same block
of $y$ (this is so because $a$ is an optimal left most alignment, see lemma \ref{AadjacentBlocks}). So between blocks of $x$ aligned with the same block
of $y$, there is at least one left out block of $x$.
 Hence the maximum proportion of blocks of $x$, which are aligned
at the same time as other blocks of $x$ to the same block of $y$, can not exceed twice the number of left out blocks
of $x$. This yields a lower bound equal to $2q_1$. This 
is as a proportion among all blocks in $x$, but we are interested
in the number as a proportion of the total number of blocks
of length $l-1$ of $x$. So, we get as lower bound $2q_1/p_{l-1}$,
where $p_{l-1}$ is the proportion of blocks of $x$ which have length
$l-1$. Adding the blocks in $x$ which are left out  and the blocks
which are aligned with several blocks of $y$, we get:
\begin{equation}
\label{trivial}
p_{l-1}^{III}\leq\frac{3q_1}{p_{l-1}}.
\end{equation}
By $D^n(\delta_0)$, we have that $p_{l-1}\geq (1/3)-\delta_0$, so that
together with \ref{trivial}, we obtain:
$$p_{l-1}^{III}\leq \frac{3q_1}{(1/3)-\delta_0}.$$
By using the above inequality in \ref{heini} we obtain:
\begin{equation}
\label{heini2}
\mathrm{E}[L(x^*,y)-L(x,y)]\geq
\frac{p_{l-1,l}+p_{l-1,l+1}}{p_{l-1,l-1}+p_{l-1,l}+p_{l-1,l+1}}
\left(1-\frac{3q_1}{(1/3)-\delta_0}  
\right).
\end{equation}
Next we are going to investigate the effect of decreasing 
 a randomly chosen block of length $l+1$ by one unit. 
The score can decrease when the
selected block of $x$ of length $l+1$
is aligned with a block of length $l+1$ of $y$. If it is aligned
with one block and that block has length $l$ or $l-1$, then there 
is no decrease. This leads to:
$$\mathrm{E}[L(\tilde{x},y)-L(x^*,y)|G_{l+1,I}]\geq 
-\frac{p_{l+1,l+1}}{p_{l+1,l-1}+p_{l+1,l}+p_{l+1,l+1}},$$
where $G_{l+1,I}$ denotes the event that the block of length
$l+1$ chosen is aligned one block with one block.
When the selected block of $x$ of length $l+1$ is aligned with
 several blocks of $y$ then the score decreases  by one unit.
 When the selected block of length $l+1$ in $x$ is left out
or is aligned at the same time as other blocks of $x$ to
the same block of $y$ then there is no decrease. This leads to:
\begin{equation}
\label{tilde*}
\mathrm{E}[L(\tilde{x},y)-L(x^*,y)]\geq 
-\mathrm{P}(G_{l+1,I})\frac{p_{l+1,l+1}}{p_{l+1,l-1}+p_{l+1,l}+p_{l+1,l+1}}
-\mathrm{P}(G_{l+1,II}),
\end{equation}
where $G_{l+1,II}$ denotes the event that the selected block 
of length $l+1$ is aligned with several blocks of $y$ at the same time.
Let $p_{l+1}$ denote the total proportion of blocks of length
$l+1$ among all blocks of $x$. Let
$p_{l+1,I}$ denote the proportion among all blocks of $x$ of length
$l+1$, of blocks which are aligned one to one. There is a proportion
of $q_1$ totally left out blocks in $x$. At most a proportion
 $\delta_0$ are at the end of the alignment a contiguous group
of left out blocks. That means, (assuming $q_1\geq \delta_0$),
the proportion of left out blocks in $x$ which are not adjacent
to another left out block of $x$ is at least $q_1-\delta_0$.
Going with each left out block which is not adjacent to another left
out block, there is at least one adjacent block which is aligned
together with several other blocks of $x$ to the same block of $y$.
This gives a lower bound for the blocks of $x$ 
which are not aligned one block to one block of
 $\delta_0+2(q_1-\delta_0)$. This is taken as a proportion among all blocks
of $x$. This gives among all blocks of length $l+1$ a proportion of
at least:
$$\frac {\delta_0+2(q_1-\delta_0)}{(1/3+\delta_0)},$$
since by the event $D^n(\delta_0)$ we know that among all blocks of $x$ the proportion
of the blocks of length $l+1$ is less than $(1/3)+\delta_0$.
 Hence,
\begin{equation}
\label{PGl+1I}
\mathrm{P}(G_{l+1,I})\leq 1-\frac {\delta_0+2(q_1-\delta_0)}{(1/3+\delta_0)}.
\end{equation}
Next let us note that we can give an upper bound for the number
of blocks of $x$ aligned with several blocks of $y$. Since we never
have several blocks aligned with several blocks, we have that
the number of blocks of $x$ aligned with several blocks of 
$y$ is not more than the total number of left out blocks
of $y$. This is so because between two  blocks 
aligned with the same block there is always at least one left
out block. The proportion of left out blocks in $y$ is $q_2$.
but this is taken as proportion among all the blocks of 
$y$. Since the total amount of blocks in $x$ and $y$
could not be exactly the same, that number can get slightly
changed when we report it as proportion
of the total number of blocks in $x$. Let $p_{l+1}$ denote the proportion among the blocks of $x$ which are of length $l+1$. We have thus that the probability to select a block of length $l+1$ of $x$ which is aligned with several blocks of $y$ is less or equal to
\begin{equation}\label{pl+1}
{\rm P}(G_{l+1,II}) \le \frac{q_2 N^Y_n}{p_{l+1} N^X_n}.
\end{equation}
By the event $D^n(\delta_0)$ we have
\begin{equation}\label{pl+1bound}
p_{l+1} \ge \frac{1}{3}-\delta_0
\end{equation}
and by the event $G^n(\delta_0)$ we have
\begin{equation}\label{eventoGn}
\frac{N^Y_n}{N^X_n} \le 1+\delta_0.
\end{equation}
Applying now \ref{pl+1bound} and \ref{eventoGn} to \ref{pl+1}, we find
\begin{equation}\label{PGl+1II}
\mathrm{P}(G_{l+1,II})\leq q_2\cdot\frac{1+\delta_0}{(1/3)-\delta_0}
\end{equation}
Finally, using inequalities \ref{PGl+1II}, \ref{PGl+1I}
in \ref{tilde*} we get:
\begin{eqnarray}
\mathrm{E}[L(\tilde{x},y)-L(x^*,y)]&\geq& -\left(1-\frac {\delta_0+2(q_1-\delta_0)}{(1/3+\delta_0)}\right)\frac{p_{l+1,l+1}}{p_{l+1,l-1}+p_{l+1,l}+p_{l+1,l+1}}\nonumber\\
&&-q_2\cdot\frac{1+\delta_0}{(1/3)-\delta_0}\label{tilde*2},
\end{eqnarray}
Using inequalities \ref{heini2} and \ref{tilde*2} together we find:
\begin{eqnarray}
\mathrm{E}[L(\tilde{x},y)-L(x,y)] &\geq& \mathrm{E}[L(\tilde{x},y)-L(x^*,y)]+\mathrm{E}[L(x^*,y)-L(x,y)] \nonumber \\
&\geq& \frac{p_{l-1,l}+p_{l-1,l+1}}{p_{l-1,l-1}+p_{l-1,l}+p_{l-1,l+1}} \left(1-\frac{3q_1}{(1/3)-\delta_0}\right) \nonumber \\
&-&\left(1-\frac {\delta_0+2(q_1-\delta_0)}{(1/3+\delta_0)}\right) \frac{p_{l+1,l+1}}{p_{l+1,l-1}+p_{l+1,l}+p_{l+1,l+1}}\nonumber\\
&&-q_2\cdot\frac{1+\delta_0}{(1/3)-\delta_0}, \label{strange}
\end{eqnarray}
Note next that we can apply lemma \ref{delta1delta2}
with $\Delta= n^{0.6}$ because of $C^n$
and $\delta_1,\delta_2\leq \delta_0$ thanks to $E^n(\delta_0)$. Hence we find that:
$$|q_1-q_2|\leq 1.5|\delta_0|+\frac{4ln^{0.6}}{n}$$
We assume that $n$ is large enough  so that:
$$|q_1-q_2|\leq 2|\delta_0|.$$
With the last inequality holding, we get from 
lemma \ref{approachOpt} that if inequality \ref{expressionI*} holds, then
\ref{Hq1+Hq2} should be satisfied. By the event  $E^n(\epsilon_2)$,
we have that \ref{Hq1+Hq2} can not be satisfied. Hence,
the inequality  \ref{expressionI*} cannot hold, which implies
that the expression on the left side of \ref{expressionI*}
is larger or equal to $\epsilon_1$. 
Together with inequality \ref{strange}, this implies
that:
\[\mathrm{E}[L(\tilde{x},y)-L(x,y)]\geq \epsilon_1.\quad\eop\]

\noindent {\bf \large Acknowledgments}
\vspace{12pt}

\noindent The authors would like to thank the support of the German Science Foundation (DFG) through the International Graduate College "Stochastics and Real World Models" (IRTG 1132) at Bielefeld University and  through the Collavorative Research Center 701 "Spectral Structures and Topological Methods in Mathematics" (CRC 701) at Bielefeld University.

\end{document}